\documentclass[10.5pt,a4paper]{article}

\usepackage{amssymb}
\usepackage{amsmath}
\usepackage{tensor}
\newtheorem{Th}{Theorem}

\newtheorem{Rm}{Remark}

\newcommand{\be}{\begin{equation}}
\newcommand{\ee}{\end{equation}}
\newcommand{\bes}{\begin{equation*}}
\newcommand{\ees}{\end{equation*}}
\newcommand{\R}{\mathbb{R}}

\newcommand\res{\mathop{\hbox{\vrule height 7pt width .5pt depth 0pt
\vrule height .5pt width 6pt depth 0pt}}\nolimits}

\newcommand{\reset}{\setcounter{equation}{0}\setcounter{Th}{0}\setcounter{Prop}{0}\setcounter{Co}{0}\setcounter{Lma}{0}\setcounter{Rm}{0}}

\def\al{\alpha}
\def\la{\lambda}
\def\eps{\varepsilon}
\def\Om{\Omega}
\def\pro{\pi_{\vec{n}}}
\def\bn{\vec{n}}

\def\bh{\vec{h}}
\def\ba{\vec{a}}
\def\bib{\vec{b}}
\def\bH{\vec{H}}
\def\bL{\vec{L}}
\def\bR{\vec{R}}

\def\bX{\vec{X}}
\def\bV{\vec{V}}

\def\bB{\vec{B}}
\def\bp{\vec{\Phi}}
\def\bT{\vec{T}}

\def\bul{\bullet}

\def\res{\mathop{\hbox{\vrule height 7pt width .5pt 
depth 0pt\vrule height .5pt width 6pt depth 0pt}}\nolimits}

\begin{document}
\title{Noether's Theorem and the Willmore Functional}
\author{Yann Bernard\footnote{Departement Mathematik, ETH-Zentrum, 8093 Z\"urich, Switzerland.}}
\date{ }
\maketitle

{\bf Abstract :} {\it Noether's theorem and the invariances of the Willmore functional are used to derive conservation laws that are satisfied by the critical points of the Willmore energy subject to generic constraints. We recover in particular previous results independently obtained by R. Capovilla and J. Guven, and by T. Rivi\`ere. Several examples are considered in details.}


\reset

\section{Introduction}

Prior to establishing herself as a leading German mathematician of the early 20$^\text{th}$ century through her seminal work in abstract algebra, Emmy Noether had already made a significant contribution to variational calculus and its applications to physics. Proved in 1915 and published in 1918 \cite{Noe}, what was to become known as {\it Noether's theorem}, is a fundamental tool in modern theoretical physics and in the calculus of variations \cite{GM, Kos, Run}. Generalizing the idea of constants of motion found in classical mechanics, Noether's theorem provides a deep connection between symmetries and conservation laws. It is a recipe to construct a divergence-free vector field from a solution of a variational problem whose corresponding action (i.e. energy) is invariant under a continuous symmetry. For example, in 1-dimensional problems where the independent  variable represents time, these vector fields are quantities which are conserved in time, such as the total energy, the linear momentum, or the angular momentum. We now precisely state one version of Noether's theorem\footnote{further generalizations may be found {\it inter alia} in \cite{Kos, Run}.}.
\medskip
Let $\Om$ be an open subset of $\mathcal{D}\subset\R^s$, and let $\mathcal{M}\subset\R^m$. Suppose that
$$
L\,:\,\Big\{(x,q,p)\,\big|\,(x,q)\in\mathcal{D}\times\mathcal{M}\;,\;p\in T_q\mathcal{M}\otimes T^*_x\mathcal{D}  \Big\}\;\longmapsto\;\R
$$
is a continuously differentiable function. Choosing a $C^1$ density measure $d\mu(x)$ on $\Om$, we can define the {\it action functional} 
$$
\mathcal{L}(u)\;:=\;\int_\Om L(x,u(x),du(x))\,d\mu(x)
$$
on the set of maps $u\in C^1(\Om,\mathcal{M})$. A tangent vector field $X$ on $\mathcal{M}$ is called an {\it infinitesimal symmetry} for $\mathcal{L}$ if it satisfies
$$
\dfrac{\partial L}{\partial q^i}(x,q,p)X^i(q)+\dfrac{\partial L}{\partial p^i_\al}(x,q,p)\dfrac{\partial X^i}{\partial q^j}(q)p^j_\al\;=\;0\:.
$$ 
\begin{Th}
Let $X$ be a Lipschitz tangent vector field on $\mathcal{M}$ which is an infinitesimal symmetry for the action $\mathcal{L}$. If $u:\Om\rightarrow\mathcal{M}$ is a critical point of $\mathcal{L}$, then
\be\label{noether}
\sum_{\al=1}^{s}\,\dfrac{\partial}{\partial x^\al}\bigg(\rho(x)X^j(u)\dfrac{\partial L}{\partial p^j_\al}(x,u,du) \bigg)\;=\;0\:,
\ee
where $\{x^\al\}_{\al=1,\ldots,s}$ are coordinates on $\Om$ such that $d\mu(x)=\rho(x)dx^1\cdot\cdot\cdot dx^s$.
\end{Th}
Equation (\ref{noether}) is the conservation law associated with the symmetry represented by $X$. The quantity
$$
\rho(x)X^j(u)\dfrac{\partial L}{\partial p^j_\al}(x,u,du)
$$
is often called {\it Noether current}, especially in the physics literature.  \\

Whether in the form given above or in analogous forms, Noether's theorem has long been recognized as a fundamental tool in variational calculus. In the context of harmonic map theory, Noether's theorem was first used by \cite{Raw} in the mid 1980s. A few years later, several authors \cite{YMC, KRS, Sha} have independently used it to replace the harmonic map equation into spheres, where derivatives of the solution appear in a quadratic way, by an equation in divergence form, where derivatives of the solution appear in a linear way. This gives a particularly helpful analytical edge when studying harmonic maps with only very weak regularity hypotheses. Fr\'ed\'eric H\'elein made significant contributions to the analysis of harmonic maps using conservation laws via Noether's theorem \cite{Hel}. In the same vein, Tristan Rivi\`ere used conservation laws to study conformally invariant variational problems \cite{Riv1}.  \\

We will in this paper also make use of Noether's theorem\footnote{not directly in the form (\ref{noether}), but the spirit behind our derivations is the same.}, this time in the context of fourth-order geometric problems in connection with the {\it Willmore functional}. We now briefly recall the main historical landmarks that led to the discovery -- and rediscovery, indeed -- of the Willmore functional. \\

Imagine that you had at your disposal the bow of a violin and a horizontal thin metallic plate covered with grains of sand. What would you observe if you were to rub the bow against the edge of the plate? In 1680, the English philosopher and scientist Robert Hooke was the first to try to answer this question (then posed in slightly different experimental terms). Some 120 years later, the German physicist and musician Ernst Chladni repeated the experiment in a systematic way \cite{Chl}. Rubbing the bow with varying frequency, he observed that the grains of sand arrange themselves in remarkable patterns -- nowadays known as {\it Chladni figures}. Those who witnessed Chladni's experiment were fascinated by the patterns, as was in 1809 the French emperor Napol\'eon I. Eager to understand the physical phenomenon at the origin of the Chladni figures, the emperor mandated Pierre-Simon de Laplace of the Acad\'emie des Sciences to organize a competition whose goal would be to provide a mathematical explanation for the figures. The winner would receive one kilogram of solid gold. Joseph-Louis Lagrange discouraged many potential candidates as he declared that the solution of the problem would require the creation of a new branch of mathematics. Only two contenders remained in the race: the very academic Sim\'eon-Denis Poisson and one autodidactic outsider: Sophie Germain. It is unfortunately impossible to give here a detailed account of the interesting events that took place in the following years (see \cite{Dah}). In 1816, Sophie Germain won the prize -- which she never claimed. Although Germain did not answer Napol\'eon's original question, and although she did not isolate the main phenomenon responsible for the Chladni figures, namely resonance, her work proved fundamental, for, as predicted by Lagrange, she laid down the foundations of a whole new branch of applied mathematics: the theory of elasticity of membranes. For the sake of brevity, one could synthesize Germain's main idea by isolating one single decisive postulate which can be inferred from her work \cite{Ger}. Having found her inspiration in the works of Daniel Bernoulli \cite{DBer} and Leonhard Euler \cite{Eul} on the elastica (flexible beams), Sophie Germain postulates that the density of elastic energy stored in a thin plate is proportional to the square of the mean curvature\footnote{Incidentally, the notion of mean curvature was first defined and used in this context; it is a creation which we owe to Germain.} $H$. In other words, the elastic energy of a bent thin plate $\Sigma$ can be expressed in the form
$$
\int_{\Sigma}H^2(p)d\sigma(p)\:,
$$
where $d\sigma$ denotes the area-element. In the literature, this energy is usually referred to as {\it Willmore energy}. It bears the name of the English mathematician Thomas Willmore who rediscovered it in the 1960s \cite{Wil1}. Prior to Willmore and after Germain, the German school of geometers led by Wilhelm Blaschke considered and studied the Willmore energy in the context of conformal geometry. Blaschke observed that minimal surfaces minimize the Willmore energy and moreover that the Willmore energy is invariant under conformal transformations of $\R^3\cup\{\infty\}$. In his nomenclature, critical points of the Willmore energy were called {\it conformal minimal surfaces} \cite{Bla}. Gerhard Thomsen, a graduate student of Blaschke, derived the Euler-Lagrange equation corresponding to the Willmore energy \cite{Tho} (this was further generalized to higher codimension in the 1970s by Joel Weiner \cite{Wei}). It is a fourth-order nonlinear partial differential equation for the immersion. Namely, let $\bp:\Sigma\rightarrow\R^{m\ge3}$ be a smooth immersion of an oriented surface $\Sigma$. The pull-back metric $g:=\bp^*g_{\R^3}$ is represented in local coordinates with components $g_{ij}$. We let $\nabla_j$ denote the corresponding covariant derivative. The second fundamental form is the normal valued 2-tensor with components $\bh_{ij}:=\nabla_i\nabla_j\bp$. Its half-trace is the mean curvature vector $\bH:=\dfrac{1}{2}\bh^{j}_{j}$. The {\it Willmore equation} reads
\be\label{will0}
\Delta_\perp\bH+\big(\bh^{i}_{j}\cdot\bH\big)\bh^{j}_{i}-2|\bH|^2\bH\;=\;\vec{0}\:,
\ee
where $\Delta_\perp$ is the negative covariant Laplacian for the connection $\nabla$ in the normal bundle derived from the ambient scalar product in $\R^m$. Note, in passing, that it is not at all clear how one could define a weak solution of (\ref{will0}) using only the requirement that $\bH$ be square-integrable (i.e. that the Willmore energy be finite). \\

The Willmore energy appears in various areas of science: general relativity, as the main contributor to the Hawking mass \cite{Haw} ; in cell biology (see below) ; in nonlinear elasticity theory \cite{FJM} ; in optical design and lens crafting \cite{KR} ; in string theory, in the guise of a string action \`a la Polyakov \cite{Pol}. As mentioned earlier, the Willmore energy also plays a distinguished role in conformal geometry, where it has given rise to many interesting problems and where it has stimulated too many elaborate works to be cited here. We content ourselves with mentioning the remarkable tours de force of Fernando Marques and Andr\'e Neves \cite{MN} to solve the celebrated Willmore conjecture stating that, up to M\"obius transformations, the Clifford torus\footnote{obtained by rotating a circle of radius 1 around an axis located at a distance $\sqrt{2}$ of its center.} minimizes the Willmore energy amongst immersed tori in $\R^3$. \\

Aiming at solving the Willmore conjecture, Leon Simon initiated the ``modern" variational study of the Willmore functional \cite{Sim} when proving the existence of an embedded torus into $\R^{m\ge3}$ minimizing the $L^2$ norm of the second fundamental form. As this norm does not provide any control of the $C^1$ norm of the surface, speaking of ``immersion" is impossible. Simon thus had to weaken the geometric notion of immersion, and did so by using varifolds and their local approximation by biharmonic graphs. In the following years, this successful ``ambient approach" was used by various authors \cite{BK, KS1, KS2, KS3} to solve important questions about Willmore surfaces. \\

Several authors \cite{CDDRR, Dal, KS2, Pal, Rus} have observed that the Willmore equation in codimension 1 is cognate with a certain divergence form. We will prove below (Theorem \ref{Th1}) a pointwise equality to that effect in any codimension. The versions found in the aforementioned works are weaker in the sense that they only identify an integral identity. \\
In 2006, Tristan Rivi\`ere \cite{Riv2} showed that the fourth-order Willmore equation (\ref{will0}) can be written in divergence form and eventually recast as a system two of second-order equations enjoying a particular structure useful to the analysis of the critical points of the Willmore energy. This observation proved to be decisive in the resolution of several questions pertaining to Willmore surfaces \cite{YBer, BR1, BR2, BR3, KMR, MR, Riv2, Riv3, Riv4}. It also led to the so-called ``parametric approach" of the problem. In contrast with the ambient approach where surfaces are viewed as subsets of $\R^m$, the parametric approach favors viewing surfaces as images of (weak) immersions, and the properties of these immersions become the analytical point of focus. \\

We briefly review the results in \cite{Riv2} (in codimension 1), adapting slightly the original notation and statements to match the orientation of our paper. With the same notation as above, the first conservation law in \cite{Riv2} states that a smooth immersion $\bp:\Sigma\rightarrow\R^3$ is a critical point of the Willmore functional if and only if
\be\label{ri1}
\nabla_j\big(\nabla^j\bH-2(\bn\cdot\nabla^j\bH)\bn+|\bH|^2\nabla^j\bp\big)\;=\;\vec{0}\:,
\ee
where $\bn$ is the outward unit normal. \\
Locally about every point, (\ref{ri1}) may be integrated to yield a function $\bL\in\R^3$ satisfying
$$
|g|^{-1/2}\epsilon^{kj}\nabla_k\bL\;=\;\nabla^j\bH-2(\bn\cdot\nabla^j\bH)\bn+|\bH|^2\nabla^j\bp\:,
$$
where $|g|$ is the volume element of the pull-back metric $g$, and $\epsilon^{kj}$ is the Levi-Civita symbol. The following equations hold:
\be\label{ri2}
\left\{\begin{array}{rcl}
\nabla_j\big(|g|^{-1/2}\epsilon^{kj}\bL\times\nabla_k\bp-\bH\times\nabla^j\bp\big)&=&\vec{0}\\[1ex]
\nabla_j\big(|g|^{-1/2}\epsilon^{kj}\bL\cdot\nabla_k\bp\big)&=&0\:.
\end{array}\right.
\ee
These two additional conservation laws give rise (locally about every point) to two potentials $\bR\in\R^3$ and $S\in\R$ satisfying
\bes
\left\{\begin{array}{rcl}
\nabla_k\bR&=&\bL\times\nabla_k\bp-|g|^{1/2}\epsilon_{kj}\bH\times\nabla^j\bp\\[1ex]
\nabla_k S&=&\bL\cdot\nabla_k\bp\:.
\end{array}\right.
\ees
A computation shows that these potentials are related to each other via the {system
\be\label{ri3}
\left\{\begin{array}{rcl}
|g|^{1/2}\Delta_g S&=&\epsilon^{jk}\partial_j\bn\cdot\partial_k\bR\\[1ex]
|g|^{1/2}\Delta_g\bR&=&\epsilon^{jk}\Big[\partial_j\bn\,\partial_kS+\partial_j\bn\times\partial_k\bR\Big]\:.
\end{array}\right.
\ee
This system is linear in $S$ and $\bR$. It enjoys the particularity of being written in flat divergence form, with the right-hand side comprising Jacobian-type terms. The Willmore energy is, up to a topological constant, the $W^{1,2}$-norm of $\bn$. For an immersion $\bp\in W^{2,2}\cap W^{1,\infty}$, one can show that $S$ and $\bR$ belong to $W^{1,2}$. Standard Wente estimates may thus be performed on (\ref{ri3}) to thwart criticality and regularity statements ensue \cite{BR1, Riv2}. Furthermore, one verifies that (\ref{ri3}) is stable under weak limiting process, which has many nontrivial consequences \cite{BR1,BR3}. \\

In 2013, the author found that the divergence form and system derived by Rivi\`ere can be obtained by applying Noether's principle to the Willmore energy\footnote{The results were first presented at Oberwolfach in July 2013 during the mini-workshop {\it The Willmore functional and the Willmore conjecture}.}. The translation, rotation, and dilation invariances of the Willmore energy yield via Noether's principle the conservations laws (\ref{ri1}) and (\ref{ri2}). 

\begin{Th}\label{Th1}
Let $\bp:\Sigma\rightarrow\R^m$ be a smooth immersion of an oriented surface $\Sigma$. Introduce the quantities 
$$
\left\{\begin{array}{lcl}
\vec{\mathcal{W}}&:=&\Delta_\perp\bH+\big(\bh^{i}_{j}\cdot\bH\big)\bh^{j}_{i}-2|\bH|^2\bH\\[1ex]
\vec{T}^j&:=&\nabla^j\bH-2\pro\nabla^j\bH+|\bH|^2\nabla^j\bp\:,
\end{array}\right.
$$
where $\pro$ denotes projection onto the normal space. \\
Via Noether's theorem, the invariance of the Willmore energy by translations, rotations, and dilations in $\R^m$ imply respectively the following three conservation laws:
\be\label{laws}
\left\{\begin{array}{rcl}
\nabla_j\vec{T}^j&=&-\,\vec{\mathcal{W}}\\[1ex]
\nabla_j\big(\bT^j\wedge\bp+\bH\wedge\nabla^j\bp\big)&=&-\,\vec{\mathcal{W}}\wedge\bp\\[1ex]
\nabla_j\big(\bT^j\cdot\bp\big)&=&-\,\vec{\mathcal{W}}\cdot\bp\:.
\end{array}\right.
\ee
In particular, the immersion $\bp$ is Willmore if and only if the following conservation holds:
\bes
\nabla_j\big(\nabla^j\bH-2\pro\nabla^j\bH+|\bH|^2\nabla^j\bp  \big)\;=\;\vec{0}\:.
\ees 
\end{Th}

For the purpose of local analysis, one may apply Hodge decompositions in order to integrate the three conservation laws (\ref{laws}). Doing so yields ``potential" functions related to each other in a very peculiar way, which we state below. The somewhat unusual notation -- the price to pay to work in higher codimension -- is clarified in Section \ref{nota}. 
\begin{Th}\label{Th2}
Let $\bp:D^2\rightarrow\R^m$ be a smooth\footnote{in practice, this strong hypothesis is reduced to $\bp\in W^{2,2}\cap W^{1,\infty}$ without modifying the result.} immersion of the flat unit disk $D^2\subset\R^2$. We denote by $\bn$ the Gauss-map, by $g:=\bp^*g_{\R^m}$ the pull-back metric, and by $\Delta_g$ the associated negative Laplace-Beltrami operator. Suppose that $\bp$ satisfies the fourth-order equation
$$
\Delta_\perp\bH+\big(\bh^{i}_{j}\cdot\bH\big)\bh^{j}_{i}-2|\bH|^2\bH\;=\;\vec{\mathcal{W}}\:,
$$
for some given $\vec{\mathcal{W}}$. Let $\bV$, $\bX$, and $Y$ solve the problems
$$
\Delta_g\bV\;=\;-\,\vec{\mathcal{W}}\qquad,\qquad \Delta_g\bX\;=\;\nabla^j\bV\wedge\nabla_j\bp\qquad,\qquad\Delta_gY\;=\;\nabla^j\bV\cdot\nabla_j\bp\:.
$$
Then $\bp$ is a solution of the second-order equation
\be\label{eqphi}
|g|^{1/2}\Delta_g\bp\;=\;-\,\epsilon^{jk}\Big[\partial_jS\partial_k\bp+\partial_j\bR\bul\partial_k\bp\Big]+|g|^{1/2}\big(\nabla^jY\nabla_j\bp+\nabla^j\bX\bul\nabla_j\bp  \big)\:,
\ee
where $S$ and $\bR$ satisfy the system 
\be\label{thesys000}
\left\{\begin{array}{rcl}
|g|^{1/2}\Delta_g S&=&\epsilon^{jk}\partial_j(\star\,\bn)\cdot\partial_k\bR
\,+\,|g|^{1/2}\nabla_j\big((\star\,\bn)\cdot\nabla^j\bX\big)\\[1ex]
|g|^{1/2}\Delta_g\bR&=&\epsilon^{jk}\Big[\partial_j(\star\,\bn)\partial_kS+\partial_j(\star\,\bn)\bul\partial_k\bR\Big]+|g|^{1/2}\nabla_j\big((\star\,\bn)\nabla^jY+(\star\,\bn)\bul\nabla^j\bX\big)\:.
\end{array}\right.
\ee
\end{Th}
In the special case when $\bp$ is Willmore, we have $\mathcal{\vec{W}}\equiv\vec{0}$, and we may choose $\bV$, $\bX$, and $Y$ to identically vanish. Then (\ref{thesys000}) becomes the conservative Willmore system found originally in \cite{Riv2}. \\

Although perhaps at first glance a little cryptic, Theorem \ref{Th2} turns out to be particularly useful for local analytical purposes. If the given $\mathcal{\vec{W}}$ is sufficiently regular, the non-Jacobian terms involving $Y$ and $\bX$ on the right-hand side of (\ref{thesys000}) form a subcritical perturbation of the Jacobian terms involving $S$ and $\bR$ (see \cite{BWW1} for details). From an analytic standpoint, one is left with studying a linear system of Jacobian-type. Wente estimates provide fine regularity information on the potential functions $S$ and $\bR$, which may, in turn, be bootstrapped into (\ref{eqphi}) to yield regularity information on the immersion $\bp$ itself. \\

In \cite{DDW}, the authors study Willmore surfaces of revolution. They use the invariances of the Willmore functional to recast the Willmore ODE in a form that is a special case of (\ref{ri1}). Applying Noether's principle to the Willmore energy had already been independently done and used in the physics community \cite{CG, Mue}. As far as the author understands, these references are largely unknown in the analysis and geometry community. One goal of this paper is to bridge the gap, as well as to present results which do not appear in print. The author hopes it will increase in the analysis/geometry community the visibility of results, which, he believes, are useful to the study of fourth-order geometric problems associated with the Willmore energy. \\
For the sake of brevity, the present work focuses only on computational derivations and on examples. A second work \cite{BWW1} written jointly with Glen Wheeler and Valentina-Mira Wheeler will shortly be available. It builds upon the reformulations given in the present paper to derive various local analytical results.  

\paragraph{Acknowledgments.} The author is grateful to Daniel Lengeler for pointing out to him \cite{CG} and \cite{Mue}. The author would also like to thank Hans-Christoph Grunau, Tristan Rivi\`ere, and Glen Wheeler for insightful discussions. The excellent working conditions of the welcoming facilities of the Forschungsinstitut f\"ur Mathematik at the ETH in Z\"urich are duly acknowledged.

\section{Main Result}

After establishing some notation in Section \ref{nota}, the contents of Theorem \ref{Th1} and of Theorem \ref{Th2} will be proved simultaneously in Section \ref{proof}.  

\subsection{Notation}\label{nota}

In the sequel, $\bp:\Sigma\rightarrow\R^{m\ge3}$ denotes a smooth immersion of an oriented surface $\Sigma$ into Euclidean space. The induced metric is $g:=\bp^*g_{\R^m}$ with components $g_{ij}$ and with volume element $|g|$. The components of the second fundamental form are denoted $\vec{h}_{ij}$. The mean curvature is $\bH:=\dfrac{1}{2}\,g^{ij}\vec{h}_{ij}$. At every point on $\Sigma$, there is an oriented basis $\{\bn_\al\}_{\al=1,\ldots,m-2}$ of the normal space. We denote by $\pro$ the projection on the space spanned by the vectors $\{\bn_\al\}$, and by $\pi_T$ the projection on the tangent space (i.e. $\pi_T+\pro=\text{id}$). The Gauss map $\bn$ is the $(m-2)$-vector defined via
$$
\star\,\bn\;:=\;\dfrac{1}{2}\,|g|^{-1/2}\epsilon^{ab}\nabla_a\bp\wedge\nabla_b\bp\:,
$$
where $\star$ is the usual Hodge-star operator, and $\epsilon^{ab}$ is the Levi-Civita symbol\footnote{Recall that the Levi-Civita is {\it not} a tensor. It satisfies $\epsilon_{ab}=\epsilon^{ab}$.} with components $\epsilon^{11}=0=\epsilon^{22}$ and $\epsilon^{12}=1=-\epsilon^{21}$. Einstein's summation convention applies throughout. We reserve the symbol $\nabla$ for the covariant derivative associated with the metric $g$. Local flat derivatives will be indicated by the symbol $\partial$. \\

\noindent
As we work in any codimension, it is helpful to distinguish scalar quantities from vector quantities. For this reason, we append an arrow to the elements of $\Lambda^p(\R^m)$, for all $p>0$. The scalar product in $\R^m$ is denoted by a dot. We also use dot to denote the natural extension of the scalar product in $\R^m$ to multivectors (see \cite{Fed}).\\ 
Two operations between multivectors are useful. The interior multiplication $\res$ maps the pair comprising a $q$-vector $\gamma$ and a $p$-vector $\beta$ to the $(q-p)$-vector $\gamma\res\beta$. It is defined via
\bes
\langle \gamma\res\beta\,,\alpha\rangle\;=\;\langle \gamma\,,\beta\wedge\alpha\rangle\:\qquad\text{for each $(q-p)$-vector $\al$.}
\ees
Let $\al$ be a $k$-vector. The first-order contraction operation $\bul$ is defined inductively through 
\bes
\al\bul\beta\;=\;\al\res\beta\:\:\qquad\text{when $\beta$ is a 1-vector}\:,
\ees
and
\bes
\al\bul(\beta\wedge\gamma)\;=\;(\al\bul\beta)\wedge\gamma\,+\,(-1)^{pq}\,(\al\bul\gamma)\wedge\beta\:,
\ees
when $\beta$ and $\gamma$ are respectively a $p$-vector and a $q$-vector.

\subsection{Variational Derivations}\label{proof}

Consider a variation of the form:
\bes
\bp_t\;:=\;\bp\,+\,t\big(A^j\nabla_j\bp+\bB\big)\:,
\ees
for some $A^j$ and some normal vector $\vec{B}$. 
We have
\bes
\nabla_i\nabla_j\bp\;=\;\bh_{ij}\:.
\ees
Denoting for notational convenience by $\delta$ the variation at $t=0$, we find:
\bes
\delta\nabla_j\bp\;\equiv\;\nabla_j\delta\bp\;=\;(\nabla_jA^s)\nabla_s\bp+A^s\bh_{js}+\nabla_j\bB\:.
\ees
Accordingly, we find
\begin{eqnarray*}
\pro\nabla^j\delta\nabla_j\bp&=&2(\nabla_jA^s)\bh^j_{s}+A^s\pro\nabla^j\bh_{js}+\pro\nabla^j\pro\nabla_j\bB+\pro\nabla^j\pi_T\nabla_j\bB\\
&=&2(\nabla_jA^s)\bh^j_{s}+2A^s\nabla_s\bH+\Delta_\perp\bB+\pro\nabla^j\pi_T\nabla_j\bB\:,
\end{eqnarray*}
where we have used the definition of the normal Laplacian $\Delta_\perp$ and the contracted Codazzi-Mainardi equation
\bes
\nabla^j\bh_{js}\;=\;2\nabla_s\bH\:.
\ees
Since $\bB$ is a normal vector, one easily verifies that
\bes
\pi_T\nabla_j\bB\;=\;-\,(\bB\cdot\bh_j^s)\,\nabla_s\bp\:,
\ees
so that
\bes
\pro\nabla^j\pi_T\nabla_j\bB\;=\;-\,(\bB\cdot\bh_j^s)\,\bh_s^j\:.
\ees
Hence, the following identity holds
\be\label{eq1}
\pro\nabla^j\delta\nabla_j\bp\;=\;2(\nabla_jA^s)\bh_{js}+2A^s\nabla_s\bH+\Delta_\perp\bB-(\bB\cdot\bh_j^s)\,\bh_s^j\:.
\ee
Note that
\be\label{eq2}
\delta g^{ij}\;=\;-2\nabla^jA^i+2\bB\cdot\bh^{ij}\:,
\ee
which, along with (\ref{eq1}) and (\ref{eq2}), then gives
\begin{eqnarray*}
\delta|\bH|^2&=&\bH\cdot\delta\nabla^j\nabla_j\bp\;\;=\;\;\bH\cdot\big[(\delta g^{ij})\partial_i\nabla_j\bp+\nabla^j\delta\nabla_j\bp\big]\\[1ex]
&=&\bH\cdot\big[(\delta g^{ij})\vec{h}_{ij}+\pro\nabla^j\delta\nabla_j\bp\big]\\[1ex]
&=&\bH\cdot\big[\Delta_\perp\bB+(\bB\cdot\bh^i_j)\bh^j_i+2A^j\nabla_j\bH  \big]\:.
\end{eqnarray*}
Finally, since
\bes
\delta|g|^{1/2}\;=\;|g|^{1/2}\big[\nabla_jA^j-2\bB\cdot\bH  \big]\:,
\ees
we obtain
\begin{eqnarray*}
\delta\big(|\bH|^2|g|^{1/2}\big)&=&|g|^{1/2}\Big[\bH\cdot\Delta_\perp\bB+(\bB\cdot\bh^i_j)\bh^j_i-2(\bB\cdot\bH)|H|^2+\nabla_j\big(|\bH|^2A^j\big)  \Big]\\[1ex]
&=&|g|^{1/2}\Big[\bB\cdot\vec{\mathcal{W}}+\nabla_j\big(\bH\cdot\nabla^j\bB-\bB\cdot\nabla^j\bH+|\bH|^2A^j\big) \Big]\:,
\end{eqnarray*}
where
\bes
\vec{\mathcal{W}}\;:=\;\Delta_\perp\bH+(\bH\cdot\bh^i_j)\bh^j_i-2|\bH|^2\bH\:.
\ees
Therefore,
\be\label{diffen}
\delta\int_{\Sigma_0}|\bH|^2\;=\;\int_{\Sigma_0}\Big[\bB\cdot\vec{\mathcal{W}}+\nabla_j\big(\bH\cdot\nabla^j\bB-\bB\cdot\nabla^j\bH+|\bH|^2A^j\big) \Big]\:.
\ee
This identity holds for every piece of surface $\Sigma_0\subset\Sigma$. We will now consider specific deformations which are known to preserve the Willmore energy (namely translations, rotations, and dilations \cite{BYC3}), and thus for which the right-hand side of (\ref{diffen}) vanishes. 

\paragraph{Translations.}

We consider a deformation of the form
\bes
\bp_t\;=\;\bp+t\vec{a}\qquad\text{for some fixed $\ba\in\R^m$}\:.
\ees
Hence
\bes
\bB\;=\;\pro\ba\qquad\text{and}\qquad A^j\;=\;\vec{a}\cdot\nabla^j\bp\:.
\ees
This gives
\begin{eqnarray*}
&&\bH\cdot\nabla^j\bB-\bB\cdot\nabla^j\bH+  |\bH|^2A^j\\[1ex]
&=&\ba\cdot\Big[(\bH\cdot\nabla^j\bn_\al-\bn_\al\cdot\nabla^j\bH)\bn_\al+H^\al\nabla^j\bn_\al+|\bH|^2\nabla^j\bp\Big]\\[1ex]
&=&\ba\cdot\Big[\nabla^j\bH-2\pro\nabla^j\bH+|\bH|^2\nabla^j\bp\Big]\:,
\end{eqnarray*}
so that (\ref{diffen}) yields
\bes
\ba\cdot\int_{\Sigma_0}\vec{\mathcal{W}}+\nabla_j\Big[\nabla^j\bH-2\pro\nabla^j\bH+|\bH|^2\nabla^j\bp  \Big]\;=\;0\:.
\ees
As this holds for all $\ba$ and all $\Sigma_0$, letting
\be\label{defT}
\bT^j\;:=\;\nabla^j\bH-2\pro\nabla^j\bH+|\bH|^2\nabla^j\bp
\ee
gives
\be\label{trans}
\nabla_j\bT^j\;=\;-\,\vec{\mathcal{W}}\:.
\ee
This is equivalent to the conservation law derived in \cite{Riv2} in the case when $\vec{\mathcal{W}}=\vec{0}$ and when the induced metric is conformal with respect to the identity. At the equilibrium, i.e. when $\vec{\mathcal{W}}$ identically vanishes, $\vec{T}^j$ plays in the problem the role of stress-energy tensor. \\

For future convenience, we formally introduce the Hodge decomposition\footnote{Naturally, this is only permitted when working locally, or on a domain whose boundary is contractible to a point.}
\be\label{defLL}
\bT^j\;=\;\nabla^j\bV+|g|^{-1/2}\epsilon^{kj}\nabla_k\bL\:,
\ee
for some $\bL$ and some $\bV$ satisfying
\be\label{defV}
-\,\Delta_g\bV\;=\;\vec{\mathcal{W}}\:.
\ee

\paragraph{Rotations.}

We consider a deformation of the form
\bes
\bp_t\;=\;\bp+t\star(\bib\wedge\bp)\qquad\text{for some fixed $\bib\in\Lambda^{m-2}(\R^m)$}\:.
\ees
In this case, we have
\bes
B^\al\;=\;-\,\bib\cdot\star(\bn_\al\wedge\bp)\qquad\text{and}\qquad A^j\;=\;-\,\bib\,\cdot\star\big(\nabla^j\bp\wedge\bp\big)\:.
\ees
Hence
\begin{eqnarray*}
&&\bH\cdot\nabla^j\bB-\bB\cdot\nabla^j\bH+  |\bH|^2A^j\\[1ex]
&&\hspace{-1cm}=\;\;-\,\bib\cdot\star\Big[(\bH\cdot\nabla^j\bn_\al-\bn_\al\cdot\nabla^j\bH)(\bn_\al\wedge\bp)+H^\al\nabla^j(\bn_\al\wedge\bp)+|\bH|^2\nabla^j\bp\wedge\bp\Big]\\
&&\hspace{-1cm}=\;\;-\,\bib\cdot\star\big(\bT^j\wedge\bp+\bH\wedge\nabla^j\bp\big)\:,
\end{eqnarray*}
where we have used the tensor $\bT^j$ defined in (\ref{defT}). 
Putting this last expression in (\ref{diffen}) and proceeding as in the previous paragraph yields the pointwise equalities
\begin{eqnarray}\label{ach1}
\vec{\mathcal{W}}\wedge\bp&=&-\,\nabla_j\big(\bT^j\wedge\bp+\bH\wedge\nabla^j\bp\big)\nonumber\\
&\stackrel{\text{(\ref{defLL})}}{\equiv}&-\,\nabla_j\big(\nabla^j\bV\wedge\bp+|g|^{-1/2}\epsilon^{kj}\nabla_k\bL\wedge\bp+\bH\wedge\nabla^j\bp\big)\nonumber\\
&=&-\,\nabla_j\big(\nabla^j\bV\wedge\bp-|g|^{-1/2}\epsilon^{kj}\bL\wedge\nabla_k\bp+\bH\wedge\nabla^j\bp\big)\nonumber\\
&=&-\Delta_g\bV\wedge\bp-\nabla^j\bV\wedge\nabla_j\bp+\nabla_j\big(|g|^{-1/2}\epsilon^{kj}\bL\wedge\nabla_k\bp-\bH\wedge\nabla^j\bp\big)\:.
\end{eqnarray}
Owing to (\ref{defV}), we thus find
\bes
\nabla_j\big(|g|^{-1/2}\epsilon^{kj}\bL\wedge\nabla_k\bp-\bH\wedge\nabla^j\bp\big)\;=\;\nabla^j\bV\wedge\nabla_j\bp\:.
\ees
It will be convenient to define two 2-vectors $\bX$ and $\bR$ satisfying the Hodge decomposition
\be\label{defR}
|g|^{-1/2}\epsilon^{kj}\bL\wedge\nabla_k\bp-\bH\wedge\nabla^j\bp\;=\;\nabla^j\bX+|g|^{-1/2}\epsilon^{kj}\nabla_k\bR\:,
\ee
with thus
\be\label{defX}
\Delta_g\bX\;=\;\nabla^j\bV\wedge\nabla_j\bp\:.
\ee

\paragraph{Dilations.}

We consider a deformation of the form
\bes
\bp_t\;=\;\bp+t\la\bp\qquad\text{for some fixed $\la\in\R$}\:,
\ees
from which we obtain
\bes
B^\al\;=\;\la\,\bn_\al\cdot\bp\qquad\text{and}\qquad A^j\;=\;\la\,\nabla^j\bp\cdot\bp\:.
\ees
Hence
\begin{eqnarray*}
&&\bH\cdot\nabla^j\bB-\bB\cdot\nabla^j\bH+  |\bH|^2A^j\\[1ex]
&&\hspace{-1cm}=\;\;\la\Big[(\bH\cdot\nabla^j\bn_\al-\bn_\al\cdot\nabla^j\bH)(\bn_\al\cdot\bp)+H^\al\nabla^j(\bn_\al\cdot\bp)+|\bH|^2\nabla^j\bp\cdot\bp  \Big]\\
&&\hspace{-1cm}=\;\;\la\,\bT^j\cdot\bp\:,
\end{eqnarray*}
where we have used that $\bH\cdot\nabla^j\bp=0$, and where $\bT^j$ is as in (\ref{defT}). \\
Putting this last expression in (\ref{diffen}) and proceeding as before gives the pointwise equalities
\begin{eqnarray}\label{ach2}
\vec{\mathcal{W}}\cdot\bp&=&-\,\nabla_j\big(\bT^j\cdot\bp\big)\\
&\equiv&-\,\nabla_j\big(\nabla^j\bV\cdot\bp+|g|^{-1/2}\epsilon^{kj}\nabla_k\bL\cdot\bp  \big)\nonumber\\
&=&-\,\Delta_g\bV\cdot\bp-\nabla^j\bV\cdot\nabla_j\bp+\nabla_j\big(|g|^{-1/2}\epsilon^{kj}\bL\cdot\nabla_j\bp\big)\nonumber\:.
\end{eqnarray}
Hence, from (\ref{defV}), we find
\bes
\nabla_j\big(|g|^{-1/2}\epsilon^{kj}\bL\cdot\nabla_k\bp\big)\;=\;\nabla^j\bV\cdot\nabla_j\bp\:.
\ees
We again use a Hodge decomposition to write
\be\label{defS}
|g|^{-1/2}\epsilon^{kj}\bL\cdot\nabla_k\bp\;=\;\nabla^jY+|g|^{-1/2}\epsilon^{kj}\nabla_k S\:,
\ee
where
\be\label{defY}
\Delta_g Y\;=\;\nabla^j\bV\cdot\nabla_j\bp\:.
\ee

\medskip

Our next task consists in relating to each other the ``potentials" $\bR$ and $S$ defined above. Although this is the fruit of a rather elementary computation, the result it yields has far-reaching consequences and which, as far as the author knows, has no direct empirical justification. 
Recall (\ref{defR}) and (\ref{defS}), namely:
\be\label{defRS}
\left\{\begin{array}{lcl}
\nabla_k\bR&=&\bL\wedge\nabla_k\bp\,-\,|g|^{1/2}\epsilon_{kj}\big(\bH\wedge\nabla^j\bp+\nabla^j\bX\big)\\[1ex]
\nabla_k S&=&\bL\cdot\nabla_k\bp\,-\,|g|^{1/2}\epsilon_{kj}\nabla^jY\:.
\end{array}\right.
\ee
Define the Gauss map 
\bes
\star\,\bn\;:=\;\dfrac{1}{2}\,|g|^{-1/2}\epsilon^{ab}\nabla_a\bp\wedge\nabla_b\bp\:.
\ees
We have
\begin{eqnarray*}
(\star\,\bn)\cdot\nabla_k\bR&=&\dfrac{1}{2}\,|g|^{-1/2}\epsilon^{ab}\big(\nabla_a\bp\wedge\nabla_b\bp\big)\cdot\Big[\bL\wedge\nabla_k\bp\,-\,|g|^{1/2}\epsilon_{kj}\big(\bH\wedge\nabla^j\bp+\nabla^j\bX\big)  \Big]\\[1ex]
&=&|g|^{-1/2}\epsilon^{ab}g_{bk}\bL\cdot\nabla_a\bp\,-\,|g|^{1/2}\epsilon_{kj}(\star\,\bn)\cdot\nabla^j\bX\\[1ex]
&=&|g|^{-1/2}\epsilon^{ab}g_{bk}\big(\nabla_aS+|g|^{1/2}\epsilon_{aj}\nabla^jY   \big)\,-\,|g|^{1/2}\epsilon_{kj}(\star\,\bn)\cdot\nabla^j\bX\\[1ex]
&=&|g|^{1/2}\epsilon_{bk}\nabla^bS+\nabla_kY\,-\,|g|^{1/2}\epsilon_{kj}(\star\,\bn)\cdot\nabla^j\bX\:,
\end{eqnarray*}
where we have used that $\bH$ is a normal vector, along with the elementary identities
\bes
|g|^{-1/2}\epsilon^{ab}g_{bk}\;=\;|g|^{1/2}\epsilon_{bk}g^{ab}\qquad\text{and}\qquad\epsilon^{ab}\epsilon_{aj}\;=\;\delta^{b}_{j}\:.
\ees
The latter implies
\be\label{nablaS}
\nabla^jS\;=\;|g|^{-1/2}\epsilon^{jk}\big((\star\,\bn)\cdot\nabla_k\bR-\nabla_kY\big)\,+\,(\star\,\bn)\cdot\nabla^j\bX\:.
\ee
Analogously, we find\footnote{$
(\omega_1\wedge\omega_2)\bul(\omega_3\wedge\omega_4)\;=\;(\omega_2\cdot\omega_4)\omega_1\wedge\omega_3-(\omega_2\cdot\omega_3)\omega_1\wedge\omega_4-(\omega_1\cdot\omega_4)\omega_2\wedge\omega_3+(\omega_1\cdot\omega_3)\omega_2\wedge\omega_4\:.
$}
\begin{eqnarray*}
(\star\,\bn)\bul\nabla_k\bR&=&\dfrac{1}{2}\,|g|^{-1/2}\epsilon^{ab}\big(\nabla_a\bp\wedge\nabla_b\bp\big)\bul\Big[\bL\wedge\nabla_k\bp\,-\,|g|^{1/2}\epsilon_{kj}\big(\bH\wedge\nabla^j\bp+\nabla^j\bX\big)  \Big]\\[1ex]
&=&|g|^{-1/2}\epsilon^{ab}g_{bk}\nabla_a\bp\wedge\bL\,+\,|g|^{-1/2}\epsilon^{ab}(\bL\cdot\nabla_a\bp)(\nabla_b\bp\wedge\nabla_k\bp)\\
&&-\:\epsilon^{ab}\epsilon_{kb}\nabla_a\bp\wedge\bH\,-\,|g|^{1/2}\epsilon_{kj}(\star\,\bn)\bul\nabla^j\bX\\[1ex]
&=&|g|^{1/2}\epsilon_{kb}\bL\wedge\nabla^b\bp\,-\,(\star\,\bn)(\bL\cdot\nabla_k\bp)\,-\,\nabla_k\bp\wedge\bH\,-\,|g|^{1/2}\epsilon_{kj}(\star\,\bn)\bul\nabla^j\bX\\[1ex]
&=&|g|^{1/2}\epsilon_{kj}\nabla^j\bR+\nabla_k\bX-(\star\,\bn)\big(\nabla_kS+|g|^{1/2}\epsilon_{kj}\nabla^jY  \big)\,-\,|g|^{1/2}\epsilon_{kj}(\star\,\bn)\bul\nabla^j\bX\:.
\end{eqnarray*}
It hence follows that there holds
\be\label{nablaR}
\nabla^j\bR\;=\;|g|^{-1/2}\epsilon^{kj}\big((\star\,\bn)\nabla_kS+(\star\,\bn)\bul\nabla_k\bR-\nabla_k\bX\big)\,+\,(\star\,\bn)\nabla^jY\,+\,(\star\,\bn)\bul\nabla^j\bX\:.
\ee
Applying divergence to each of (\ref{nablaS}) and (\ref{nablaR}) gives rise to the {\it conservative Willmore system}:
\be\label{conswillsys}
\left\{\begin{array}{lcl}
|g|^{1/2}\Delta_g S&=&\epsilon^{jk}\partial_j(\star\,\bn)\cdot\partial_k\bR
\,+\,|g|^{1/2}\nabla_j\big((\star\,\bn)\cdot\nabla^j\bX\big)\\[1ex]
|g|^{1/2}\Delta_g\bR&=&\epsilon^{jk}\Big[\partial_j(\star\,\bn)\partial_kS+\partial_j(\star\,\bn)\bul\partial_k\bR\Big]+|g|^{1/2}\nabla_j\big((\star\,\bn)\nabla^jY+(\star\,\bn)\bul\nabla^j\bX\big)\:,
\end{array}\right.
\ee
where $\partial_j$ denotes the derivative in flat local coordinates. \\

\noindent
This system is to be supplemented with (\ref{defX}) and (\ref{defY}), which, in turn, are solely determined by the value of the Willmore operator $\vec{\mathcal{W}}$ via equation (\ref{defV}). There is furthermore another useful equation to add to the conservative Willmore system, namely one relating the potentials $S$ and $\bR$ back to the immersion $\bp$. We now derive this identity. Using (\ref{defRS}), it easily follows that
\begin{eqnarray*}
&&\hspace{-2cm}\epsilon^{km}\big(\nabla_k\bR+|g|^{1/2}\epsilon_{kj}\nabla^j\bX\big)\bul\nabla_m\bp\\[1ex]
&=&\epsilon^{km}\big(\bL\wedge\nabla_k\bp-|g|^{1/2}\epsilon_{kj}\bH\wedge\nabla^j\bp\big)\bul\nabla_m\bp\\[1ex]
&=&\epsilon^{km}\Big[\big(\bL\cdot\nabla_m\bp\big)\nabla_k\bp\,-\,g_{mk}\bL\,-\,|g|^{1/2}\epsilon_{km}\bH\Big]\\[1ex]
&=&\epsilon^{km}\big(\nabla_mS+|g|^{1/2}\epsilon_{mj}\nabla^jY\big)\nabla_k\bp\,-\,2\,|g|^{1/2}\bH\:.
\end{eqnarray*}
Since $\Delta_g\bp=2\bH$, we thus find
\be\label{law4}
\partial_j\Big[\epsilon^{jk}\big(S\partial_k\bp+\bR\bul\partial_k\bp   \big)+|g|^{1/2}\nabla^j\bp\Big]\;=\;|g|^{1/2}\big(\nabla^jY\nabla_j\bp+\nabla^j\bX\bul\nabla_j\bp  \big)\:.
\ee
At the equilibrium (i.e. when $\vec{\mathcal{W}}=\vec{0}$), the right-hand side of the latter identically vanishes and we recover a conservation law. With the help of a somewhat tedious computation, one verifies that this conservation law follows from the invariance of the Willmore energy by inversion and Noether's theorem. To do so, one may for example consider an infinitesimal variation of the type $\delta\bp\;=\;|\bp|^2\ba-2(\bp\cdot\ba)\bp$, for some fixed constant vector $\ba\in\R^m$.\\

We summarize our results in the following pair of systems. 
\be\label{side}
\left\{\begin{array}{rcl}\vec{\mathcal{W}}&=&\Delta_\perp\bH+(\bH\cdot\vec{h}_j^i)\vec{h}^j_i-2|\bH|^2\bH\\[1ex]
\Delta_g\vec{V}&=&-\,\vec{\mathcal{W}}\\[1ex]
\Delta_g\vec{X}&=&\nabla^j\bV\wedge\nabla_j\bp\\[1ex]
\Delta_g{Y}&=&\nabla^j\bV\cdot\nabla_j\bp
\end{array}\right.
\ee
and
\be\label{thesys}
\left\{\begin{array}{rcl}
|g|^{1/2}\Delta_g S&=&\epsilon^{jk}\partial_j(\star\,\bn)\cdot\partial_k\bR
\,+\,|g|^{1/2}\nabla_j\big((\star\,\bn)\cdot\nabla^j\bX\big)\\[1ex]
|g|^{1/2}\Delta_g\bR&=&\epsilon^{jk}\Big[\partial_j(\star\,\bn)\partial_kS+\partial_j(\star\,\bn)\bul\partial_k\bR\Big]+|g|^{1/2}\nabla_j\big((\star\,\bn)\nabla^jY+(\star\,\bn)\bul\nabla^j\bX\big)\\[1ex]
|g|^{1/2}\Delta_g\bp&=&-\,\epsilon^{jk}\Big[\partial_jS\partial_k\bp+\partial_j\bR\bul\partial_k\bp\Big]+|g|^{1/2}\big(\nabla^jY\nabla_j\bp+\nabla^j\bX\bul\nabla_j\bp  \big)\:.
\end{array}\right.
\ee

\noindent
In the next section, we will examine more precisely the structure of this system through several examples. 

\medskip
\begin{Rm}
Owing to the identities
\bes
\vec{u}\bul\vec{v}\;=\;(\star\vec{u})\times\vec{v}\quad\text{and}\quad \vec{u}\bul\vec{w}\;=\;\star\big[(\star\vec{u})\times(\star\vec{w})\big]\:\quad\text{for}\:\:\:\vec{u}\in\Lambda^2(\R^3), \vec{v}\in\Lambda^1(\R^3), \vec{w}\in\Lambda^2(\R^3)\:,
\ees
we can in $\R^3$ recast the above systems as
\be\label{side1}
\left\{\begin{array}{rcl}\vec{\mathcal{W}}&=&\Delta_\perp\bH+(\bH\cdot\vec{h}_j^i)\vec{h}^j_i-2|\bH|^2\bH\\[1ex]
\Delta_g\vec{V}&=&-\,\vec{\mathcal{W}}\\[1ex]
\Delta_g\vec{X}&=&\nabla^j\bV\times\nabla_j\bp\\[1ex]
\Delta_g{Y}&=&\nabla^j\bV\cdot\nabla_j\bp
\end{array}\right.
\ee
and
\be\label{thesys1}
\left\{\begin{array}{rcl}
|g|^{1/2}\Delta_g S&=&\epsilon^{jk}\partial_j\bn\cdot\partial_k\bR
\,+\,|g|^{1/2}\nabla_j\big(\bn\cdot\nabla^j\bX\big)\\[1ex]
|g|^{1/2}\Delta_g\bR&=&\epsilon^{jk}\Big[\partial_j\bn\,\partial_kS+\partial_j\bn\times\partial_k\bR\Big]+|g|^{1/2}\nabla_j\big(\bn\,\nabla^jY+\bn\times\nabla^j\bX\big)\\[1ex]
|g|^{1/2}\Delta_g\bp&=&-\,\epsilon^{jk}\Big[\partial_jS\,\partial_k\bp+\partial_j\bR\times\partial_k\bp\Big]+|g|^{1/2}\big(\nabla^jY\nabla_j\bp+\nabla^j\bX\times\nabla_j\bp  \big)\:.
\end{array}\right.
\ee
In this setting, $\bX$ and $\bR$ are no longer 2-vectors, but rather simply vectors of $\R^3$.
\end{Rm}

\section{Examples}

\subsection{Willmore Immersions}

A smooth immersion $\bp:\Sigma\rightarrow\R^{m\ge3}$ of an oriented surface $\Sigma$ with induced metric $g=\bp^*g_{\R^m}$ and corresponding mean curvature vector $\bH$, is said to be Willmore if it is a critical point of the Willmore energy $\int_\Sigma|\bH|^2d\text{vol}_g$. They are known \cite{Wil2,Wei} to satisfy the Euler-Lagrange equation
$$
\Delta_\perp\bH+(\bH\cdot\vec{h}_j^i)\vec{h}^j_i-2|\bH|^2\bH\;=\;\vec{0}\:.
$$
In the notation of the previous section, this corresponds to the case $\vec{\mathcal{W}}=\vec{0}$. According to (\ref{side}), we have the freedom to set $\vec{V}$, $\vec{X}$, and $Y$ to be identically zero. The Willmore equation then yields the second-order system in divergence form
\be\label{thesyswillmore}
\left\{\begin{array}{rcl}
|g|^{1/2}\Delta_g S&=&\epsilon^{jk}\partial_j(\star\,\bn)\cdot\partial_k\bR
\\[1ex]
|g|^{1/2}\Delta_g\bR&=&\epsilon^{jk}\Big[\partial_j(\star\,\bn)\partial_kS+\partial_j(\star\,\bn)\bul\partial_k\bR\Big]\\[1ex]
|g|^{1/2}\Delta_g\bp&=&-\,\epsilon^{jk}\Big[\partial_jS\partial_k\bp+\partial_j\bR\bul\partial_k\bp\Big]\:.
\end{array}\right.
\ee
This system was originally derived by Rivi\`ere in \cite{Riv2}. For notational reasons, the detailed computations were carried out only in local conformal coordinates, that is when $g_{ij}=\text{e}^{2\la}\delta_{ij}$, for some conformal parameter $\la$. The analytical advantages of the Willmore system (\ref{thesyswillmore}) have been exploited in numerous works \cite{BR1, BR2, BR3, Riv2, Riv3, Riv4}. The flat divergence form of the operator $|g|^{1/2}\Delta_g$ and the Jacobian-type structure of the right-hand side enable using fine Wente-type estimates in order to produce non-trivial local information about Willmore immersions (see aforementioned works). 

\begin{Rm}
Any Willmore immersion will satisfy the system (\ref{thesyswillmore}). The converse is however not true. Indeed, in order to derive (\ref{thesyswillmore}), we first obtained the existence of some ``potential" $\bL$ satisfying the first-order equation
\be\label{ceteq}
|g|^{-1/2}\epsilon^{kj}\nabla_k\bL\;=\;\nabla^j\bH-2\pro\nabla^j\bH+|\bH|^2\nabla^j\bp\:.
\ee
In doing so, we have gone from the Willmore equation, which is second-order for $\bH$, to the above equation, which is only first-order in $\bH$, thereby introducing in the problem an extraneous degree of freedom. As we shall see in the next section, (\ref{ceteq}) is in fact equivalent to the conformally-constrained Willmore equation, which, as one might suspect, is the Willmore equation supplemented with an additional degree of freedom appearing in the guise of a Lagrange multiplier.  
\end{Rm}

\subsection{Conformally-Constrained Willmore Immersions}

Varying the Willmore energy $\int_\Sigma|\bH|^2d\text{vol}_g$ in a fixed conformal class (i.e. with infinitesimal, smooth, compactly supported, conformal variations) gives rise to a more general class of surfaces called {\it conformally-constrained Willmore surfaces} whose corresponding Euler-Lagrange equation \cite{BPP, KL, Sch} is expressed as follows. 
Let $\vec{h}_0$ denote the trace-free part of the second fundamental form, namely
\bes
\vec{h}_0\;:=\;\vec{h}\,-\,\bH g\:.
\ees
A conformally-constrained Willmore immersion $\bp$ satisfies
\be\label{cwe}
\Delta_\perp\bH+(\bH\cdot\vec{h}_j^i)\vec{h}^j_i-2|\bH|^2\bH\;=\;\big(\vec{h}_0\big)_{ij}q^{ij}\:,
\ee
where $q$ is a transverse\footnote{i.e. $q$ is divergence-free: $\nabla^j q_{ji}=0\:\:\forall i$.}  traceless symmetric 2-form. This tensor $q$ plays the role of Lagrange multiplier in the constrained variational problem. \\

In \cite{KS3}, it is shown that under a suitable ``small enough energy" assumption, a minimizer of the Willmore energy in a fixed conformal class exists and is smooth. The existence of a minimizer without any restriction on the energy is also obtained in \cite{Riv3} where it is shown that the minimizer is either smooth (with the possible exclusion of finitely many branch points if the energy is large enough to grant their formation) or else isothermic\footnote{the reader will find in \cite {Riv5} an interesting discussion on isothermic immersions.}. One learns in \cite{Sch} that non-degenerate critical points of the Willmore energy constrained to a fixed conformal class are solutions of the conformally constrained Willmore equation. 
Continuing along the lines of \cite{Riv3}, further developments are given in \cite{Riv6}, where the author shows that if either the genus of the surface satisfies $g\le2$, or else if the Teichm\"uller class of the immersion is not hyperelliptic\footnote{A class in the Teichm\"uller space is said to be {\it hyperelliptic} if the tensor products of holomorphic 1-forms do not generate the vector space of holomorphic quadratic forms.}, then any critical point $\bp$ of the Willmore energy for $C^1$ perturbations included in a submanifold of the Teichm\"uller space is in fact analytic and (\ref{cwe}) is satisfied for some transverse traceless symmetric 2-tensor $q$. \\

The notion of conformally constrained Willmore surfaces clearly generalizes that of Willmore surfaces, obtained via all smooth compactly supported infinitesimal variations (setting $q\equiv0$ in (\ref{cwe})). In \cite{KL}, it is shown that CMC Clifford tori are conformally constrained Willmore surfaces. In \cite{BR1}, the conformally constrained Willmore equation (\ref{cwe}) arises as that satisfied by the limit of a Palais-Smale sequence for the Willmore functional. \\

Minimal surfaces are examples of Willmore surfaces, while parallel mean curvature surfaces\footnote{parallel mean curvature surfaces satisfy $\pro d\bH\equiv\vec{0}$. They generalize to higher codimension the notion of constant mean curvature surfaces defined in $\R^3$. See [YBer].} are examples of conformally-constrained Willmore surfaces\footnote{{\it a non}-minimal parallel mean curvature surface is however {\it not} Willmore (unless of course it is the conformal transform of a Willmore surface ; e.g. the round sphere).}. Not only is the Willmore energy invariant under reparametrization of the domain, but more remarkably, it is invariant under conformal transformations of $\R^m\cup\{\infty\}$.
Hence, the image of a [conformally-constrained] Willmore immersion through a conformal transformation is again a [conformally-constrained] Willmore immersion. It comes thus as no surprise  that the class of Willmore immersions [resp. conformally-constrained Willmore immersions] is considerably larger than that of immersions whose mean curvature vanishes [resp. is parallel], which is {\it not} preserved through conformal diffeomorphisms. \\

Comparing (\ref{cwe}) to the first equation in (\ref{side}), we see that $\vec{\mathcal{W}}=-\,\big(\vec{h}_0\big)_{ij}q^{ij}$. 
Because $q$ is traceless and transverse, we have
$$
\big(\vec{h}_0\big)_{ij}q^{ij}\;\equiv\;\vec{h}_{ij}q^{ij}-\vec{H}g_{ij}q^{ij}\;=\;\vec{h}_{ij}q^{ij}\;=\;\nabla_j(q^{ij}\nabla_i\bp)\:.
$$
Accordingly, choosing $\nabla^j\vec{V}=-\,q^{ij}\nabla_i\bp$ will indeed solve the second equation (\ref{side}). Observe next that
$$
\nabla^j\bV\cdot\nabla_j\bp\;=\;-\,q^{ij}g_{ij}\;=\;0\:,
$$
since $q$ is traceless. Putting this into the fourth equation of (\ref{side}) shows that we may choose $Y\equiv0$. Furthermore, as $q$ is symmetric, it holds
$$
\nabla^j\bV\wedge\nabla_j\bp\;=\;-\,q^{ij}\nabla_i\bp\wedge\nabla_j\bp\;=\;\vec{0}\:,
$$
so that the third equation in (\ref{side}) enables us to choose $\bX\equiv\vec{0}$. \\
Altogether, we see that a conformally-constrained Willmore immersion, just like a ``plain" Willmore immersion (i.e. with $q\equiv0$) satisfies the system (\ref{thesyswillmore}). In fact, it was shown in \cite{BR1} that to any smooth solution $\bp$ of (\ref{thesyswillmore}), there corresponds a transverse, traceless, symmetric 2-form $q$ satisfying (\ref{cwe}).

\subsection{Bilayer Models}
Erythrocytes (also called red blood cells) are the body's principal mean of transporting vital oxygen to the organs and tissues. The cytoplasm (i.e. the ``inside") of an erythrocyte is rich in hemoglobin, which chemically tends to bind to oxygen molecules and retain them. To maximize the possible intake of oxygen by each cell, erythrocytes  -- unlike all the other types of cells of the human body -- have no nuclei\footnote{this is true for all mammals, not just for humans. However, the red blood cells of birds, fish, and reptiles do contain a nucleus.}. The membrane of an erythrocyte is a bilayer made of amphiphilic molecules. Each molecule is made of a ``head" (rather large) with a proclivity for water, and of a ``tail" (rather thin) with a tendency to avoid water molecules. When such amphiphilic molecules congregate, they naturally arrange themselves in a bilayer, whereby the tails are isolated from water by being sandwiched between two rows of heads. The membrane can then close over itself to form a vesicle. Despite the great biochemical complexity of erythrocytes, some phenomena may be described and explained with the sole help of physical mechanisms\footnote{The most celebrated such phenomenon was first observed by the Polish pathologist Tadeusz Browicz in the late 19$^{\text{th}}$ century \cite{Bro}. Using a microscope, he noted that the luminous intensity reflected by a red blood cell varies erratically along its surface, thereby giving the impression of flickering. In the 1930s, the Dutch physicist Frits Zernicke invented the phase-constrast microscope which revealed that the erratic flickering found by Browicz is in fact the result of very minute and very rapid movements of the cell's membrane. These movements were finally explained in 1975 by French physicists Fran\c{c}oise Brochard and Jean-Fran\c{c}ois Lennon: they are the result of a spontaneous thermic agitation of the membrane, which occurs independently of any particular biological activity \cite{BL}.}. For example, to understand the various shapes an erythrocyte might assume, it is sensible to model the red blood cell by a drop of hemoglobin (``vesicle") whose membrane is made of a lipid bilayer \cite{Lip}. Such ``simple" objects, called {\it liposomes}, do exist in Nature, and they can be engineered artificially\footnote{The adjective ``simple" is to be understood with care: stable vesicles with non-trivial topology can be engineered and observed. See \cite{MB} and \cite{Sei}. }. The membrane of a liposome may be seen as a viscous fluid separated from water by two layers of molecules. Unlike a solid, it can undergo shear stress. Experimental results have however shown that vesicles are very resistant and the stress required to deform a liposome to its breaking-point is so great that in practice, vesicles evolving freely in water are never submitted to such destructive forces. One may thus assume that the membrane of a liposome is incompressible: its area remains constant. The volume it encloses also stays constant, for the inside and the outside of the vesicle are assumed to be isolated. As no shearing and no stretching are possible, one may wonder which forces dictate the shape of a liposome. To understand the morphology of liposomes, Canham \cite{Can}, Evans \cite{Eva}, and Helfrich \cite{Hef} made the postulate that, just as any other elastic material does, the membrane of a liposome tends to minimize its elastic energy. As we recalled in the introduction, the elastic energy of a surface is directly proportional to the Willmore energy. Accordingly, to understand the shape of a liposome, one would seek minimizers of the Willmore energy subject to constraints on the area and on the enclosed volume. A third constraint could be taken into account. As the area of the inner layer of the membrane is slightly smaller than the area of its outer layer, it takes fewer molecules to cover the former as it does to cover the latter. This difference, called {\it asymmetric area difference}, is fixed once and for all when the liposome forms. Indeed, no molecule can move from one layer to the other, for that would require its hydrophilic head to be, for some time, facing away from water. From a theoretical point of view, the relevant energy to study is thus:
\be\label{canham}
E\;:=\;\int_{\Sigma}H^2d\text{vol}_g+\al A(\Sigma)+\beta V(\Sigma)+\gamma M(\Sigma)\:,
\ee
where $H$ denotes the mean curvature scalar\footnote{in this section, we will content ourselves with working in codimension 1.}, $A(\Sigma)$, $V(\Sigma)$, $M(\Sigma)$ denote respectively the area, volume, and asymmetric area difference of the membrane $\Sigma$, and where $\al$, $\beta$, and $\gamma$ are three Lagrange multipliers. 
Depending on the authors' background, the energy (\ref{canham}) bears the names Canham, Helfrich, Canham-Helfrich, bilayer coupling model, spontaneous curvature model, among others.\\
This energy -- in the above form or analogous ones -- appears prominently in the applied sciences literature. It would be impossible to list here a comprehensive account of relevant references. The interested reader will find many more details in \cite{Sei} and \cite{Voi} and the references therein. In the more ``analytical" literature, the energy (\ref{canham}) is seldom found (except, of course, in the case when all three Lagrange multipliers vanish). We will, in time, recall precise instances in which it has been studied. But prior to doing so, it is interesting to pause for a moment and better understand a term which might be confusing to the mathematician reader, namely the asymmetric area difference $M(\Sigma)$. In geometric terms, it is simply the total curvature of $\Sigma$:
\be\label{defM}
M(\Sigma)\;:=\;\int_{\Sigma}H\,d\text{vol}_g\:.
\ee
This follows from the fact that the infinitesimal variation of the area is the mean curvature, and thus the area difference between two nearby surfaces is the first moment of curvature. Hence, we find the equivalent expression
\be\label{canam}
E\;=\;\int_{\Sigma}\bigg(H+\dfrac{\gamma}{2}\bigg)^{\!2}d\text{vol}_g+\bigg(\al-\dfrac{\gamma^2}{4}\bigg) A(\Sigma)+\beta V(\Sigma)\:.
\ee
This form of the bilayer energy is used, inter alia, in \cite{BWW2, Whe1}, where the constant $-\gamma/2$ is called {\it spontaneous curvature}. \\

From a purely mathematical point of view, one may study the energy (\ref{canham}) not just for embedded surfaces, but more generally for immersions. An appropriate definition for the volume $V$ must be assigned to such an immersion $\bp$. As is shown in \cite{MW}, letting as usual $\bn$ denote the outward unit-normal vector to the surface, one defines
\bes
V(\Sigma)\;:=\;\int_\Sigma\bp^*(d\mathcal{H}^3)\;=\;\int_{\Sigma}\bp\cdot\bn\,d\text{vol}_g\:,
\ees
where $d\mathcal{H}^3$ is the Hausdorff measure in $\R^3$. Introducing the latter and (\ref{defM}) into (\ref{canham}) yields
\be\label{veneer}
E\;=\;\int_{\Sigma}\big(H^2+\gamma H+\beta\,\bp\cdot\bn+\al\big)\,d\text{vol}_g\:.
\ee

We next vary the energy $E$ along a normal variation of the form $\delta\bp=\bB\equiv B\vec{n}$. Using the computations from the previous section, it is not difficult to see that
\be\label{varr1}
\delta\int_\Sigma d\text{vol}_g\;=\;-\,2\int_{\Sigma}\bB\cdot\vec{H} \,d\text{vol}_g\qquad\text{and}\qquad \delta\int_\Sigma H\,d\text{vol}_g\;=\;\int_{\Sigma}
(\bB\cdot\bn)\bigg(\dfrac{1}{2}h^i_jh^j_i-2H^2\bigg)\,d\text{vol}_g\:.
\ee
With a bit more effort, in \cite{MW}, it is shown that
\be\label{varr2}
\delta\int_\Sigma \bp\cdot\bn\, d\text{vol}_g\;=\;-\int_{\Sigma}\bB\cdot\bn\,d\text{vol}_g\:.
\ee
Putting (\ref{diffen}), (\ref{varr1}), and (\ref{varr2}) into (\ref{veneer}) yields the corresponding Euler-Lagrange equation
\be\label{ELhelf}
\vec{\mathcal{W}}\;:=\;\Delta_\perp\bH+(\bH\cdot\vec{h}_j^i)\vec{h}^j_i-2|\bH|^2\bH\;=\;2\,\al\bH+\beta\,\bn-\gamma\bigg(\dfrac{1}{2}h^i_jh^j_i-2H^2  \bigg)\,\bn\:.
\ee
We now seek a solution to the second equation in (\ref{side1}), namely a vector $\bV$ satisfying $\Delta_g\bV=-\vec{\mathcal{W}}$. To do so, it suffices to observe that
\bes
2\bH\;=\;\Delta_g\bp\qquad\text{and}\qquad \big(h^i_jh^j_i-4H^2\big)\bn\;=\;\nabla_j\Big[(h^{ij}-2Hg^{ij})\nabla_i\bp\Big]\:,
\ees
where we have used the Codazzi-Mainardi identity. 
Furthermore, it holds
\bes 
2\,\bn\;=\;|g|^{-1/2}\epsilon^{ij}\nabla_i\bp\times\nabla_j\bp\;=\;\nabla_j\big[-|g|^{-1/2}\epsilon^{ij}\bp\times\nabla_i\bp   \big]\:.
\ees
Accordingly, we may choose
\bes
\nabla^j\bV\;=\;-\,\al\,\nabla^j\bp\,+\,\dfrac{\beta}{2}\,|g|^{-1/2}\epsilon^{ij}\,\bp\times\nabla_i\bp\,+\,\dfrac{\gamma}{2}\big(h^{ij}-2Hg^{ij}\big)\nabla_i\bp\:.
\ees
Introduced into the third equation of (\ref{side1}), the latter yields
\begin{eqnarray}
\Delta_g\bX&=&\nabla^j\bV\times\nabla_j\bp\;\;=\;\;\dfrac{\beta}{2}\,|g|^{-1/2}\epsilon^{ij}\big(\bp\times\nabla_i\bp\big)\times\nabla_j\bp\nonumber\\[0ex]
&=&\dfrac{\beta}{2}\,|g|^{-1/2}\epsilon^{ij}\big[(\bp\cdot\nabla_j\bp)\nabla_i\bp-g_{ij}\bp\big]\;\;=\;\;\dfrac{\beta}{4}\,|g|^{-1/2}\epsilon^{ij}\,\nabla_j|\bp|^2\,\nabla_i\bp\:,\nonumber
\end{eqnarray}
so that we may choose
\bes\label{XHel}
\nabla^j\bX\;=\;\dfrac{\beta}{4}\,|g|^{-1/2}\epsilon^{ij}|\bp|^2\,\nabla_i\bp\:.
\ees
Then we find
\be\label{Xprop}
\bn\times\nabla^j\bX\;=\;\dfrac{\beta}{4}\,|\bp|^2\nabla^j\bp\qquad\text{and}\qquad\nabla^j\bX\times\nabla^j\bp\;=\;\dfrac{\beta}{2}\,|\bp|^2\,\bn\:.
\ee
Analogously, the fourth equation of (\ref{side1}) gives
\begin{eqnarray}\label{YHel}
\Delta_gY&=&\nabla^j\bV\cdot\nabla_j\bp\;\;=\;\;-\,2\al\,-\,\gamma H\,+\,\dfrac{\beta}{2}\,|g|^{-1/2}\epsilon^{ij}\big(\bp\times\nabla_i\bp\big)\cdot\nabla_j\bp\nonumber\\[0ex]
&=&-\,2\al\,-\,\gamma H\,+\,\beta\,\bp\cdot\bn\:.
\end{eqnarray}
With (\ref{Xprop}) and (\ref{YHel}), the system (\ref{thesys1}) becomes
\be\label{thesyshel}
\left\{\begin{array}{rcl}
|g|^{1/2}\Delta_g S&=&\epsilon^{jk}\partial_j\bn\cdot\partial_k\bR\\[1ex]
|g|^{1/2}\Delta_g\bR&=&\epsilon^{jk}\Big[\partial_j\bn\,\partial_kS+\partial_j\bn\times\partial_k\bR\Big]+|g|^{1/2}\nabla_j\bigg(\bn\,\nabla^jY+\dfrac{\beta}{4}\,|\bp|^2\nabla^j\bp \bigg)\\[1.5ex]
|g|^{1/2}\Delta_g\bp&=&-\,\epsilon^{jk}\Big[\partial_jS\,\partial_k\bp+\partial_j\bR\times\partial_k\bp\Big]+|g|^{1/2}\bigg(\nabla_j\bp\nabla^jY+\dfrac{\beta}{2}\,|\bp|^2\,\bn  \bigg)\:.
\end{array}\right.
\ee

Under suitable regularity hypotheses (e.g. that the immersion be locally Lipschitz and lie in the Sobolev space $W^{2,2}$), one can show that the non-Jacobian term in the second equation, namely
$$
|g|^{1/2}\nabla_j\bigg(\bn\,\nabla^jY+\dfrac{\beta}{4}\,|\bp|^2\nabla^j\bp \bigg)\:,
$$
is a subcritical perturbation of the Jacobian term. Analyzing (\ref{thesyshel}) becomes then very similar to analyzing the Willmore system (\ref{thesyswillmore}). Details may be found in \cite{BWW1}.  \\

In some cases, our so-far local considerations can yield global information. If the vesicle we consider has the topology of a sphere, every loop on it is contractible to a point. The Hodge decompositions which we have performed in Section II to deduce the existence of $\bV$, and subsequently that of $\bX$ and $Y$ hold globally. Integrating (\ref{YHel}) over the whole surface $\Sigma$ then gives the {\it balancing condition}:
$$
2\al A(\Sigma)+\gamma M(\Sigma)\;=\;\beta V(\Sigma)\:.
$$
This condition is well-known in the physics literature \cite{CG, Sei}.

\begin{Rm}
Another instance in which minimizing the energy (\ref{veneer}) arises is the isoperimetric problem \cite{KMR, Scy}, which consists in minimizing the Willmore energy under the constraint that the dimensionless isoperimetric ratio $\sigma:=A^3/V^2$ be a given constant in $(0,1]$. As both the Willmore energy and the constraint are invariant under dilation, one might fix the volume $V=1$, forcing the area to satisfy $A=\sigma^{1/3}$. This problem is thus equivalent to minimizing the energy (\ref{veneer}) with $\gamma=0$ (no constraint imposed on the total curvature, but the volume and area are prescribed separately). One is again led to the system (\ref{thesyshel}) and local analytical information may be inferred. 
\end{Rm}

\subsection{Chen's Problem}

An isometric immersion $\bp:N^{n}\rightarrow\R^{m>n}$ of an $n$-dimensional Riemannian manifold $N^n$ into Euclidean space is called {\it biharmonic} if the corresponding mean-curvature vector $\bH$ satisfies 
\be\label{chen}
\Delta_g\bH\;=\;\vec{0}\:.
\ee
The study of biharmonic submanifolds was initiated by B.-Y. Chen \cite{BYC1} in the mid 1980s as he was seeking a classification of the finite-type submanifolds in Euclidean spaces. Independently, G.Y. Jiang \cite{Jia} also studied (\ref{chen}) in the context of the variational analysis of the biharmonic energy in the sense of Eells and Lemaire. Chen conjectures that a biharmonic immersion is necessarily minimal\footnote{The conjecture as originally stated is rather analytically vague: no particular hypotheses on the regularity of the immersion are {\it a priori} imposed. Many authors consider only smooth immersions.}.
Smooth solutions of (\ref{chen}) are known to be minimal for $n=1$ \cite{Dim1}, for $(n,m)=(2,3)$ \cite{Dim2}, and for $(n,m)=(3,4)$ \cite{HV}. A growth condition allows a simple PDE argument to work in great generality \cite{Whe2}. Chen's conjecture has also been solved under a variety of hypotheses (see the recent survey paper \cite{BYC2}). The statement remains nevertheless open in general, and in particular for immersed surfaces in $\R^m$. In this section, we show how our reformulation of the Willmore equation may be used to recast the fourth-order equation (\ref{chen}) in a second-order system with interesting analytical features. \\

Let us begin by inspecting the tangential part of (\ref{chen}), namely,
\begin{eqnarray}\label{tgtchen}
\vec{0}&=&\pi_T\Delta_g\bH\;=\;\big(\nabla_k\bp\cdot\nabla_j\nabla^j\bH)\nabla^k\bp\nonumber\\[1ex]
&=&\Big[\nabla_j\big(\nabla_k\bp\cdot\nabla^j\bH  \big)-\vec{h}_{jk}\cdot\nabla^j\bH \Big]\nabla^k\bp\;\;=\;\;-\,\Big[\bH\cdot\nabla_j\vec{h}^{j}_{k}+2\,\vec{h}_{jk}\cdot\nabla^j\bH \Big]\nabla^k\bp\nonumber\\[1ex]
&=&-\,\Big[\nabla_k|\bH|^2+2\,\vec{h}_{jk}\cdot\nabla^j\bH \Big]\nabla^k\bp\:,
\end{eqnarray}
where we have used that $\bH$ is a normal vector, as well as the Codazzi-Mainardi identity. With the help of this equation, one obtains a decisive identity, which we now derive, and which too makes use of the Codazzi-Mainardi equation. 
\begin{eqnarray}\label{deci}
\nabla_j\Big[|\bH|^2\nabla^j\bp-2(\bH\cdot\vec{h}^{jk})\nabla_k\bp\Big]&=&-\,\nabla_j|\bH|^2\nabla^j\bp+2|\bH|^2\bH-2(\vec{h}^{jk}\cdot\nabla_j\bH)\nabla_k\bp-2(\bH\cdot\vec{h}^{jk})\vec{h}_{jk}\nonumber\\[1ex]
&\stackrel{\text{(\ref{tgtchen})}}{=}&2|\bH|^2\bH-2(\bH\cdot\vec{h}^{jk})\vec{h}_{jk}\:.
\end{eqnarray}
Note that, in general, there holds
\begin{eqnarray}\label{bidule}
\pi_T\nabla^j\bH&=&\big[\nabla_k\bp\cdot\nabla^j\bH\big]\nabla^k\bp\;\;=\;\;-\,\big[\bH\cdot\vec{h}^{j}_{k}\big]\nabla^k\bp\:.
\end{eqnarray}
An immersion whose mean curvature satisfies (\ref{chen}) has thus the property that
\begin{eqnarray*}
\Delta_\perp\bH&:=&\pro\nabla_j\pro\nabla^j\bH\;\;=\;\;\pro\Delta_g\bH\,-\,\pro\nabla_j\pi_T\nabla^j\bH\;\;=\;\;\pro\nabla_j\Big[\big(\bH\cdot\vec{h}^{j}_{k}\big)\nabla^k\bp \Big]\nonumber\\[0ex]
&=&\big(\bH\cdot\vec{h}^{j}_{k}\big)\vec{h}^{k}_{j}\:.
\end{eqnarray*}
Putting the latter into the first equation of (\ref{side}) yields
\be\label{ELchen}
\vec{\mathcal{W}}\;:=\;\Delta_\perp\bH+(\bH\cdot\vec{h}_j^i)\vec{h}^j_i-2|\bH|^2\bH\;=\;2(\bH\cdot\vec{h}^{jk})\vec{h}_{jk}-2|\bH|^2\bH\:.
\ee
We now seek a solution to the second equation in (\ref{side1}), namely a vector $\bV$ satisfying $\Delta_g\bV=-\vec{\mathcal{W}}$. To do so, it suffices to compare (\ref{ELchen}) and (\ref{deci}) to see that we may choose
\bes
\nabla^j\bV\;=\;|\bH|^2\nabla^j\bp-2(\bH\cdot\vec{h}^{jk})\nabla_k\bp\:.
\ees
Introduced into the third equation of (\ref{side1}), the latter yields immediately
\bes
\Delta_g\bX\;=\;\nabla^j\bV\wedge\nabla_j\bp\;=\;\vec{0}\:,
\ees
thereby prompting us to choosing $\bX\equiv\vec{0}$. The fourth equation of (\ref{side1}) gives
\bes
\Delta_gY\;=\;\nabla^j\bV\cdot\nabla_j\bp\;=\;-\,2|\bH|^2\:.
\ees
On the other hand, owing to (\ref{chen}) and (\ref{bidule}), one finds
\bes
\Delta_g(\bp\cdot\bH)\;=\;-\,2|\bH|^2\:,
\ees
so that we may choose $Y=\bp\cdot\bH$. Introducing these newly found facts for $\bX$ and $Y$ into the system (\ref{thesys1}) finally gives
\be\label{thesyschen}
\left\{\begin{array}{rcl}
|g|^{1/2}\Delta_g S&=&\epsilon^{jk}\partial_j(\star\,\bn)\cdot\partial_k\bR\\[1ex]
|g|^{1/2}\Delta_g\bR&=&\epsilon^{jk}\Big[\partial_j(\star\,\bn)\partial_kS+\partial_j(\star\,\bn)\bul\partial_k\bR\Big]+|g|^{1/2}\nabla_j\big((\star\,\bn)\nabla^j(\bp\cdot\bH)\big)\\[1ex]
|g|^{1/2}\Delta_g\bp&=&-\,\epsilon^{jk}\Big[\partial_jS\partial_k\bp+\partial_j\bR\bul\partial_k\bp\Big]+|g|^{1/2}\big(\nabla^j(\bp\cdot\bH)\nabla_j\bp \big)\:.
\end{array}\right.
\ee

As previously noted, the reformulation of (\ref{trans}) in the form (\ref{thesys1}) is mostly useful to reduce the nonlinearities present in (\ref{trans}). Moreover, a local analysis of (\ref{thesys1}) is possible when the non-Jacobian terms on the right-hand side are subcritical perturbations of the Jacobian terms. This is not {\it a priori} the case for (\ref{thesyschen}), and the equation (\ref{chen}) is linear to begin with. Nevertheless, the system (\ref{thesyschen}) has enough suppleness\footnote{owing mostly to the fact that the function $Y:=\bp\cdot\bH$ satisfies $\Delta_gY\le 0$ and $\Delta_g^2Y\le0$.} and enough structural features to deduce interesting analytical facts about solutions of (\ref{chen}), under mild regularity requirements. This is discussed in detail in a work to appear \cite{BWW1}\footnote{see also \cite{Whe3} which contains interesting estimates for equation (\ref{chen}).}.

\subsection{Point-Singularities}

As was shown in \cite{YBer, BR2, Riv2}, the Jacobian-type system (\ref{side}) is particularly suited to the local analysis of point-singularities. The goal of this section is not to present a detailed account of the local analysis of point-singularities -- this is one of the topics of \cite{BWW1} -- but rather to give the reader a few pertinent key arguments on how this could be done. \\

Let $\bp:D^2\setminus\{0\}\rightarrow\R^m$ be a smooth immersion of the unit-disk, continuous at the origin (the origin will be the point-singularity in question). In order to make sense of the Willmore energy of the immersion $\bp$, we suppose that $\int_{D^2}|\bH|^2d\text{vol}_g<\infty$. Our immersion is assumed to satisfy the problem
\be\label{ptprob}
\Delta_\perp\bH+(\bH\cdot\vec{h}_j^i)\vec{h}^j_i-2|\bH|^2\bH\;=\;\vec{\mathcal{W}}\qquad\text{on}\:\:D^2\setminus\{0\}\:,
\ee
where the vector $\vec{\mathcal{W}}$ is given. It may depend only on geometric quantities (as is the case in the Willmore problem or in Chen's problem), but it may also involve ``exterior" quantities (as is the case in the conformally-constrained Willmore problem). To simplify the presentation, we will not in this paper discuss the integrability assumptions that must be imposed on $\vec{\mathcal{W}}$ to carry out the procedure that will be outlined. The interested reader is invited to consult \cite{BWW1} for more details on this topic.  \\

\noindent
As we have shown in (\ref{defT}), equation (\ref{ptprob}) may be rephrased as 
\bes
\partial_j\big(|g|^{1/2}\bT^j\big)\;=\;-\,\vec{\mathcal{W}}\qquad\text{on}\:\:D^2\setminus\{0\}\:,
\ees
for some suitable tensor $\bT^j$ defined solely in geometric terms. Consider next the problem
\bes
\Delta_g\bV\;=\;-\,\vec{\mathcal{W}}\qquad\text{on}\:\:D^2\:.
\ees
As long as $\vec{\mathcal{W}}$ is not too wildly behaved, this equation will have at least one solution. Let next $\mathcal{L}_g$ satisfy
\bes
\partial_j\big(|g|^{1/2}\nabla^j\mathcal{L}_g\big)\;=\;\delta_0\qquad\text{on}\:\:D^2\:.
\ees
If the immersion is correctly chosen (e.g. $\bp\in W^{2,2}\cap W^{1,\infty}$), the solution $\mathcal{L}_g$ exists and has suitable analytical properties (see \cite{BWW1} for details). \\
We have
\be\label{poinc1}
\partial_j\Big[|g|^{1/2}\big(\bT^j-\nabla^j\bV-\vec{\beta}\,\nabla^j\mathcal{L}_g\big)\Big]\;=\;\vec{0}\qquad\text{on}\:\:D^2\setminus\{0\}\:,
\ee
for any constant $\vec{\beta}\in\R^m$, and in particular for the unique $\vec{\beta}$ fulfilling the circulation condition that
\be\label{poinc2}
\int_{\partial D^2}\vec{\nu}\cdot\big(\bT^j-\nabla^j\bV-\vec{\beta}\,\nabla^j\mathcal{L}_g\big)\;=\;0\:,
\ee
where $\vec{\nu}\in\R^2$ denotes the outer unit normal vector to the boundary of the unit-disk. This vector $\vec{\beta}$ will be called {\it residue}.\\
Bringing together (\ref{poinc1}) and (\ref{poinc2}) and calling upon the Poincar\'e lemma, one infers the existence of an element $\bL$ satisfying
\be
\bT^j-\nabla^j\bV-\vec{\beta}\,\nabla^j\mathcal{L}_g\;=\;|g|^{-1/2}\epsilon^{kj}\nabla_k\bL\:,
\ee
with the same notation as before. We are now in the position of repeating {\it mutatis mutandis} the computations derived in the previous section, taking into account the presence of the residue. We define $\bX$ and $Y$ via:
\be\label{ptside}
\left\{\begin{array}{rcl}
\Delta_g\vec{X}&=&\nabla^j\big(\bV+\vec{\beta}\mathcal{L}_g\big)\wedge\nabla_j\bp\\[1ex]
\Delta_g{Y}&=&\nabla^j\big(\bV+\vec{\beta}\mathcal{L}_g\big)\cdot\nabla_j\bp
\end{array}\right.\qquad\text{on}\:\:D^2\:.
\ee
One verifies that the following equations hold
\bes
\left\{\begin{array}{rcl}
\nabla^k\Big[\bL\wedge\nabla_k\bp\,-\,|g|^{1/2}\epsilon_{kj}\big(\bH\wedge\nabla^j\bp+\nabla^j\bX\big)\Big]&=&\vec{0}\\[2ex]
\nabla^k\Big[\bL\cdot\nabla_k\bp\,-\,|g|^{1/2}\epsilon_{kj}\nabla^jY\Big]&=&0
\end{array}\right.\qquad\text{on}\:\:D^2\setminus\{0\}\:.
\ees
Imposing suitable hypotheses on the integrability of $\vec{\mathcal{W}}$ yields that the bracketed quantities in the latter are square-integrable. With the help of a classical result of Laurent Schwartz \cite{Scw}, the equations may be extended without modification to the whole unit-disk. As before, this grants the existence of two potential functions $S$ and $\bR$ which satisfy (\ref{defRS}) and the system (\ref{thesys}) on $D^2$. The Jacobian-type/divergence-type structure of the system sets the stage for a local analysis argument, which eventually yields a local expansion of the immersion $\bp$ around the point-singularity. This expansion involves the residue $\vec{\beta}$. The procedure was carried out in details for Willmore immersions in \cite{BR2}\footnote{An equivalent notion of residue was also identified in \cite{KS2}.}, and for conformally constrained Willmore immersions in \cite{YBer}. Further considerations can be found in \cite{BWW1}.


\bigskip


\begin{thebibliography}{99}
\bibitem[BK]{BK} Bauer, M.; Kuwert, E. ``Existence of minimizing Willmore surfaces of prescribed genus." Int. Math. Res. Not. (2003), no. 10, 553--576.
\bibitem[YBer]{YBer} Bernard, Y. ``Analysis of constraint Willmore surfaces." arXiv:1211.4455.
\bibitem[BR1]{BR1} Bernard, Y.; Rivi\`ere, T. ``Local Palais-Smale sequences for the Willmore functional." Comm. Anal. Geom. 19 (2011), 563--599.
\bibitem[BR2]{BR2} Bernard, Y.; Rivi\`ere, T. ``Singularity removability at branch points for Willmore surfaces." Pacific J. Math. 265 (2013), 257--311.
\bibitem[BR3]{BR3} Bernard, Y.; Rivi\`ere, T. ``Energy quantization for Willmore surfaces and applications." Ann. of Math. 180 (2014), 87--136. 
\bibitem[BWW1]{BWW1} Bernard, Y.; Wheeler, G.; Wheeler, V.-M.; to appear. 
\bibitem[BWW2]{BWW2} Bernard, Y.; Wheeler, G.; Wheeler, V.-M.``Spherocytosis and the Helfrich model." to appear. 
\bibitem[DBer]{DBer} Bernoulli, D. 26$^{\text{th}}$ letter to Euler in ``Correspondance math\'ematique et physique de quelques c\'el\`ebres g\'eom\`etres du XVIII\`eme si\`ecle" vol. 2, Nicolaus Fuss. 
\bibitem[Bla]{Bla} Blaschke, W. ``Vorlesungen \"uber Differentialgeometrie und geometrische Grundlagen von Einsteins Relativit\"atstheorie.Ó Die Grundlehren der mathematischen Wissenschaften in Einzeldarstellungen. Bd. XXIX Differentialgeometrie der Kreise und Kugeln, bearbeitet von Gerhard Thomsen (1929).
\bibitem[BPP]{BPP} Bohle, C.; Peters, G.P.; Pinkall, U. ``Constrained Willmore surfaces." Calc. Var. Partial Differential Equations 32 (2008), 263--277.
\bibitem[BL]{BL} Brochard, F.; Lennon J.-F.``Frequency spectrum of the flicker phenomenon in erythrocytes." J. Phys. France 36 (1975), 1035--1047.
\bibitem[Bro]{Bro} Browicz, T. ``Further observation of motion phenomena on red blood cells in pathological states." Zbl. Med. Wiss. 28 (1890), no. 1, 625--627.
\bibitem[Can]{Can} Canham, P.B. ``The minimum energy of bending as a possible explanation of the biconcave shape of the human red blood cell." J. Theor. Biol. 26 (1970), 61--81.
\bibitem[CG]{CG} Capovilla, R.; Guven, J. ``Stresses in lipid membranes". J. Phys. A: Math. Gen. 35 (2002), 6233--6247.
\bibitem[BYC1]{BYC1} Chen, B.-Y. ``Some open problems and conjectures on submanifolds of finite type." Soochow J. Math. 17 (1991), 169--188.
\bibitem[BYC2]{BYC2} Chen, B.-Y. ``Recent developments of biharmonic conjecture and modified biharmonic conjectures." arXiv:1307.0245.
\bibitem[BYC3]{BYC3} Chen, B.-Y. ``Some conformal invariants of submanifolds and their applications."  Boll. Un. Mat. Ital. (4) 10 (1974), 380--385.
\bibitem[YMC]{YMC} Chen, Y.M. ``The weak solutions to the evolution problems of harmonic maps." Math. Zeit. 201 (1989), 69--74. 
\bibitem[Chl]{Chl} Chladni, E. ``Die Akustik." Leipzig (1802). 
\bibitem[CDDRR]{CDDRR} Clarenz, U.; Diewald, U.; Dziuk, G.; Rumpf, M.; Rusu, R. ``A finite element method for surface restoration with smooth boundary conditions." Comput. Aided Geom. Design 21 (2004), no. 5, 427--445.
\bibitem[Dah]{Dah} Dahan-Dalm\'edico, A. ``M\'ecanique et th\'eorie des surfaces : les travaux de Sophie Germain." Hist. Math. 14 (1987), 347--365. 
\bibitem[Dal]{Dal} Dall'Acqua, A. ``Uniqueness for the homogeneous Dirichlet Willmore boundary value problem." Ann. Glob. Anal. Geom. 42 (2012), no. 3, 411--420.
\bibitem[DDW]{DDW} Dall'Acqua, A.; Deckelnick, K.; Wheeler, G. ``Unstable Willmore surfaces of revolution subject to natural boundary conditions." Calc. Var. 48 (2013), 293--313.
\bibitem[Dim1]{Dim1} Dimitri\'c, I. ``Submanifolds of $\mathbb{E}^n$ with harmonic mean curvature vector." Bull. Inst. Math. Acad. Sinica 20 (1992), 53--65.
\bibitem[Dim2]{Dim2} Dimitri\'c, I. ``Quadric representation and submanifolds of finite type." Ph.D. Thesis, Michigan State University (1989).
\bibitem[Eul]{Eul} Euler, L. ``Methodus inveniendi lineas curvas maximi minimive proprietate gaudentes, sive solutio problematis isoperimetrici lattissimo sensu accepti." Additamentum 1 (1744). 
\bibitem[Eva]{Eva} Evans, E.A. ``Bending resistance and chemically induced moments in membrane bilayers." Biophys. J. 14 (1974), 923--931.
\bibitem[Fed]{Fed} Federer, H. ``Geometric Measure Theory." Die Grundlehren der math. Wissensch. 153. Springer (1969). 
\bibitem[FJM]{FJM} Friesecke, G.; James, R.D.; M\"uller, S. ``A theorem on geometric rigidity and the derivation of nonlinear plate theory from three-dimensional elasticity." Comm. Pure Appl. Math. 55 (2002), no. 11, 1461--1506.
\bibitem[Ger]{Ger} Germain, S. ``Recherches sur la th\'eorie des surfaces \'elastiques." Courcier (1821).
\bibitem[GM]{GM} Gotay, M.J.; Marsden, J.E. ``Stress-energy tensors and the Belinfante-Rosenfeld formula." Contemp. Math. 132 (1992), 367--392.
\bibitem[HV]{HV} Hasanis, T.; Vlachos, T. ``Hypersurfaces in $\mathbb{E}^4$ with harmonic mean curvature vector field." Math. Nachr. 172 (1995), 145--169.
\bibitem[Haw]{Haw} Hawking, S.W. ``Gravitational radiation in an expanding universe."  J. Math. Phys. 9 (1968), 598--604.
\bibitem[Hel]{Hel} H\'elein, F. ``Harmonic Maps, Conservation Laws, and Moving Frames." Cambridge Tracts in Mathematics, 150. Cambridge University Press (2002).
\bibitem[Hef]{Hef} Helfrich, W. ``Elastic properties of lipid bilayers: theory and possible experiments."  Z. Naturforsch., C28 (1973), 693--703. 
\bibitem[Jia]{Jia} Jiang, G.Y. ``2-harmonic isometric immersions between Riemannian manifolds." Chinese Ann. Math. Ser. A, 7 (1986), 130--144.
\bibitem[KR]{KR} Katzman, D.; Rubinstein, J. ``Method for the design of multifocal optical elements." U.S. Patent no. US006302540B1 (October 2001). 
\bibitem[KRS]{KRS} Keller, J.; Rubinstein, J.; Sternberg, P. ``Reaction-diffusion processes and evolution to harmonic maps." SIAM J. Appl. Math. 49 (1989), no. 6, 1722--1733.
\bibitem[KMR]{KMR} Keller, L.G.A.; Mondino, A.; Rivi\`ere, T. ``Embedded surfaces of arbitrary genus minimizing the Willmore energy under isoperimetric constraint." Arch. Rat. Mech. Anal. 212 (2014), 645--682.
\bibitem[Kos]{Kos} Kosmann-Schwarzbach, Y. ``The Noether Theorems." Sources and Studies in the History of Mathematics and Physical Sciences, Springer (2011).
\bibitem[KL]{KL} Kuwert, E.; Lorenz, J. ``On the stability of the CMC Clifford tori as constrained Willmore surfaces." Ann. Glob. Anal. Geom. 44 (2013), 23--42.
\bibitem[KS1]{KS1} Kuwert, E.; Sch\"atzle, R. ``The Willmore flow with small initial energy." J. Diff. Geom. 57 (2001), 409--441.
\bibitem[KS2]{KS2} Kuwert, E.; Sch\"atzle R. ``Removability of point singularities of Willmore surfaces." Ann. of Math. 160 (2004), no. 1, 315--357.
\bibitem[KS3]{KS3} Kuwert, E.; Sch\"atzle, R. ``Minimizers of the Willmore energy under fixed conformal class." J. Diff. Geom. 93 (2013), 471--530.
\bibitem[Lip]{Lip} Lipowsky, R. ``The conformation of membranes." Nature 345 (1991), 475--481.
\bibitem[MN]{MN} Marques, F.C.; Neves, A. ``Min-max theory and the Willmore conjecture." Ann. of Math. 179 (2014), 683--782.
\bibitem[MW]{MW} McCoy, J.; Wheeler, G. ``A classification theorem for Helfrich surfaces." Math. Annalen 357 (2013), 1485--1508. 
\bibitem[MB]{MB} Michalet, X.; Bensimon, D. ``Vesicles of toroidal topology: observed morphology and shape transformations." J. Phys. II France 5 (1995), 263--287.  
\bibitem[MR]{MR} Mondino, A.; Rivi\`ere, T. ``Willmore spheres in compact Riemannian manifolds." Adv. Math. 232 (2013), 608--676.
\bibitem[Mue]{Mue} M\"uller, M.M. ``Theoretical studies of fluid membrane mechanics."
Diss. zur Erlangung des Grades Doktor der Naturwiss. am Fachbereich Physik, Mathematik und Informatik, Mainz (2007).
\bibitem[Noe]{Noe} Noether, E. ``Invariante Variationsprobleme." Nachr. d. K\"onig. Gesellsch. d. Wiss. zu G\"ottingen, Math-Phys. Klasse (1918), 235--257.
\bibitem[Olv]{Olv} Olver, P. ``Applications of Lie Groups to Differential Equations."  Graduate Texts in Mathematics 107, Springer (1993). 
\bibitem[Pal]{Pal} Palmer, B. ``Uniqueness theorems for Willmore surfaces with fixed and free boundaries." Indiana Univ. Math. J. 49 (2000), no. 4, 1581--1601.
\bibitem[Pol]{Pol} Polyakov, A.M. ``Fine structure of strings." Nucl. Phys. B 268 (1986), no. 2, 406--412.
\bibitem[Raw]{Raw} Rawnsley, J.H. ``Noether's theorem for harmonic maps." Diff. Geom. Methods in Math. Phys. (1984), 197--202.
\bibitem[Riv1]{Riv1} Rivi\`ere, T. ``Conservation laws for conformally invariant variational problems." Invent. Math. 168 (2006), no. 1, 1--22.
\bibitem[Riv2]{Riv2} Rivi\`ere, T. ``Analysis aspects of the Willmore functional." Invent. Math. 174 (2008), no. 1, 1--45.
\bibitem[Riv3]{Riv3} Rivi\`ere, T.  ``Variational Principles for immersed Surfaces with $L^2$-bounded Second Fundamental Form.'' J. reine angew. Math. (2013).
\bibitem[Riv4]{Riv4} Rivi\`ere, T. ``Lipschitz conformal immersions from degenerating Riemann surfaces with $L^2$-bounded second fundamental forms." Adv. Calc. Var. 6 (2013), 1--31. 
\bibitem[Riv5]{Riv5} Rivi\`ere, T. ``Sequences of smooth global isothermic immersions." Comm. PDE 38 (2012), no. 2, 276--303.
\bibitem[Riv6]{Riv6} Rivi\`ere, T. ``Critical weak immersed surfaces within submanifolds of the Teichm\"uller space." arXiv:1307.5406.
\bibitem[Run]{Run} Rund, H. ``The Hamilton-Jacobi Theory in the Calculus of Variations, its Role in Mathematics and Physics." Krieger (1973).
\bibitem[Rus]{Rus} Rusu, R.E. ``An algorithm for the elastic flow of surfaces." Interfaces and Free Boundaries 7 (2005), 229--239.
\bibitem[Sch]{Sch} Sch\"atzle, R. ``Conformally constrained Willmore immersions." Adv. Calc. Var. 6 (2013), 375--390.
\bibitem[Scy]{Scy} Schygulla, J. ``Willmore minimizers with prescribed isoperimetric ratio", Arch. Rat. Mech. Anal. 203 (2012), no. 3, 901--941.
\bibitem[Scw]{Scw} Schwartz, L. ``Th\'eorie des distributions" vol 1. Hermann (1951).
\bibitem[Sei]{Sei} Seifert, U. ``Configurations of fluid membranes and vesicles." Adv. in Phys. 46 (1997), no. 1, 13--137.
\bibitem[Sha]{Sha} Shatah, J. ``Weak solutions and developments of singularities of the $SU(2)$ $\sigma$-model." Comm. Pure Appl. Math. 41 (1988), 459--469.
\bibitem[Sim]{Sim} Simon, L. ``Existence of surfaces minimizing the Willmore functional." Comm. Anal. Geom. 1 (1993), no. 2, 281--326.
\bibitem[Tho]{Tho} Thomsen, G. ``\"Uber konforme Geometrie I ; Grundlagen der konformen Flachentheorie." Ab. Math. Sem. Univ. Hamburg 3 (1924), no. 1, 31--56.
\bibitem[Voi]{Voi} Voinova, M. ``Geometrical methods in the theory of lipid membranes and cells shapes." Living State Physics, Chalmers University of Technology (2006), 1--191.
\bibitem[Wei]{Wei} Weiner, J. ``On a problem of Chen, Willmore, et al.'' Indiana Univ. Math. J. 27 (1978), no. 1, 19--35.
\bibitem[Whe1]{Whe1} Wheeler, G. ``Global analysis of the generalized Helfrich flow of closed curves immersed in $\R^n$." arXiv:1205.5939.
\bibitem[Whe2]{Whe2} Wheeler, G. ``Chen's conjecture and $\eps$-superbiharmonic sub manifolds of Riemannian manifolds." Int. J. Math. 24 (2013), no. 4, 135--141. 
\bibitem[Whe3]{Whe3} Wheeler, G. ``Gap phenomena for a class of fourth-order geometric differential operators on surfaces with boundary." Proc. AMS (to appear). arXiv:1302.4165.
\bibitem[Wil1]{Wil1} Willmore, T. J. ``Note on embedded surfaces." Ann. Stiint. Univ. ``Al. I. Cuza" Iasi. Sect. I a Mat. (N.S.) 11B (1965), 493--496.
\bibitem[Wil2]{Wil2} Willmore, T. J. ``Riemannian Geometry." Oxford University Press (1997).
\end{thebibliography}
\end{document}